\documentclass[11pt,leqno,twoside]{article}

\usepackage{amsmath,amsfonts,amssymb,latexsym,oldgerm,array,epsfig,graphicx,llabel}

\usepackage[english]{babel}

\newtheorem{Theorem}{Theorem}[section]

\newtheorem{Lemma}[Theorem]{Lemma}

\newtheorem{Proposition}[Theorem]{Proposition}

\renewcommand{\i}{\iota}

\newcommand{\be}{\begin{equation}}

\newcommand{\ee}{\end{equation}}

\newcommand{\ber}{\begin{eqnarray}}

\newcommand{\eer}{\end{eqnarray}}

\newcommand{\bew}{\begin{eqnarray*}}

\newcommand{\eew}{\end{eqnarray*}}

\newcommand{\real}{{\mathbb R}}

\newcommand{\nat}{{\Bbb N}}

\renewcommand{\.}{\hspace{-2ex} .}

\renewcommand{\,}{\mbox{\hspace{1mm}}}

\newenvironment{proof}{\vspace{3pt}\indent
                       \textsc{\em Proof:}\quad }
                        {\hfill$\square$\vspace{3pt}}

\newcommand{\seq}{\longrightarrow}

\newcommand{\const}{\mbox{\rm const}}
\usepackage[arrow,matrix]{xy}

\usepackage{amssymb,german}

\usepackage{amsbsy,amsmath}

\newtheorem{cor}[Theorem]{Corollary}

\newtheorem{definition}[Theorem]{Definition}

\newtheorem{rem}[Theorem]{Remark}

\newtheorem{example}[Theorem]{Example}

\newtheorem{construction}[Theorem]{Construction}

\newcommand{\inj}{\hookrightarrow}

\newcommand{\imp}{\Rightarrow}

\newcommand{\bs}{\boldsymbol}

\newcommand{\codim}{\mbox{codim}}

\newcommand{\hol}{\mbox{\rm Hol}}

\newcommand{\diff}{\mbox{\rm Diff}}

\newcommand{\sat}{\mbox{sat}}

\renewcommand{\ker}{\mbox{ker}}
\newcommand{\bb}{\mathbb}

\newenvironment{defn}{\begin{definition}\rm}{\end{definition}}

\newenvironment{ex}[1]{\begin{example}\rm{\bf #1}}{\hfill$\bigtriangleup$\end{example}}

\renewcommand{\labelenumi}{{(\alph{enumi})}}

\renewcommand{\labelenumii}{(\roman{enumii})}

\begin{document}

\title{\bf On graphical foliations and the global existence of Euler's multiplier}

\author{
\begin{tabular}{cc}Marco K\"uhnel\footnote{The author gratefully acknowledges support by the DFG-priority program ``Global Methods in Complex Geometry''.}& Michael Neudert\\ 
{\tiny Department of Mathematics} &{\tiny Department of Mathematics}\\
{\tiny Otto-von-Guericke-University Magdeburg} &{\tiny University of Bayreuth}\\ 
{\tiny P.O. Box 4120} &{\tiny D-95440 Bayreuth}\\
{\tiny D-39016 Magdeburg} &\\
\end{tabular}}

\maketitle

\begin{abstract}
We present a 
criterion for the global existence of
Euler's multiplier for an integrable one-form taking into account 
the corresponding 
codim-1-foliation. 
In particular, the impact of inseparable leaves is considered.
Here, we suppose that the foliation can be reduced to a graph.
The properties of this graph are crucial for the 
global existence of the
Euler's multiplier. 
As applications we investigate some special cases in which the graph
turns out to look very simple.
\end{abstract}

%

%

%

%

%


%

%

%

\section{Introduction}

%

First order differential equations of the form
\begin{equation}\label{dgl}
g(x,y)dx+h(x,y)dy=0
\end{equation}
sometimes can be transformed
to exact differential equations and be solved by finding a so called integrating factor,
i.e.\ some function $\lambda=\lambda(x,y)\not= 0$ such that there is some function $f=f(x,y)$
with the property
\[
df\ =\ \lambda\ g \ dx+\lambda\ h\ dy.
\]
In general, $\lambda$ is supposed to vanish nowhere or, at least, in a set 
with no interior points (see \cite{marik}).
Then $f$ is a first integral of equation (\ref{dgl}) and the solutions are implicitly given by
the equation
$
f(x,y)=\const.
$
It is well known that, if $|g(x_0,y_0)|+|h(x_0,y_0)| > 0$, there 
are some neighbourhood
$U$ of $(x_0,y_0)$ and nontrivial functions $\lambda, f: U\to\real$ as above.
\\[0.6ex]
In higher dimensions, the situation is more complicated.
Let $r\in{\bb N}$, $r\ge 1$, and $M$ be a real $C^r$ manifold of dimension $n+1$ and $\omega$ a $C^r$ 
one-form on $M$. In this article we assume a manifold to be second countable, i.e. it has a countable topological base. This implies that it is paracompact, hence
metrizable, and can be covered by a countable number of charts. 
Now we ask for the existence of functions $f\in C^r(M,\real)$, $\lambda\in C^{r-1}(M,\real)$  
such that $\lambda(x)\not=0$  and
\[
df=\lambda \omega
\]
in $M$.
In this case, $\lambda$ is called Euler's multiplier.
From $d(df)=0$ follows immediately that
\begin{equation}\label{notw}
\omega \land d\omega\ =\ 0
\end{equation}
is a necessary condition for the existence of such functions $\lambda$ and $f$. In this
case $\omega$ is called an integrable one-form.
From Frobenius' theorem one can derive the {\it local} existence of $\lambda$ and $f$ 
in a neighbourhood of $p\in M$ if condition (\ref{notw}) is satisfied and
$\omega(p)\not=0$ (see e.g.\ \cite{vonwestenholz}).
But even if these conditions are satisfied everywhere in $M$ and $M$ is simply connected, 
global existence is,
in general, not guaranteed. Consider the following example: Let
$M=\real^3$, represented in cylinder coordinates $(\varrho,\phi,z)\in\ [0,\infty)\times [0,2\pi)\times \real$, $\ w_\varrho:=\varrho^2$,
\[
{w_\phi:=\begin{cases} 0,& 0 \le \varrho < \frac12\\ 4(\varrho-\frac12)^2(\varrho-1), 
& \frac12 \le \varrho < 1\\ \varrho-1, & 1\le \varrho\end{cases}}\quad\quad
{w_z:=\begin{cases} (\frac14-\varrho^2)^2,& 0 \le \varrho < \frac12\\ 0, & \frac12 \le 
\varrho \end{cases}}
\]
%
%
%
\begin{figure}[h!]
\begin{minipage}[h]{7cm}
\begin{center}
\includegraphics[width=24ex]{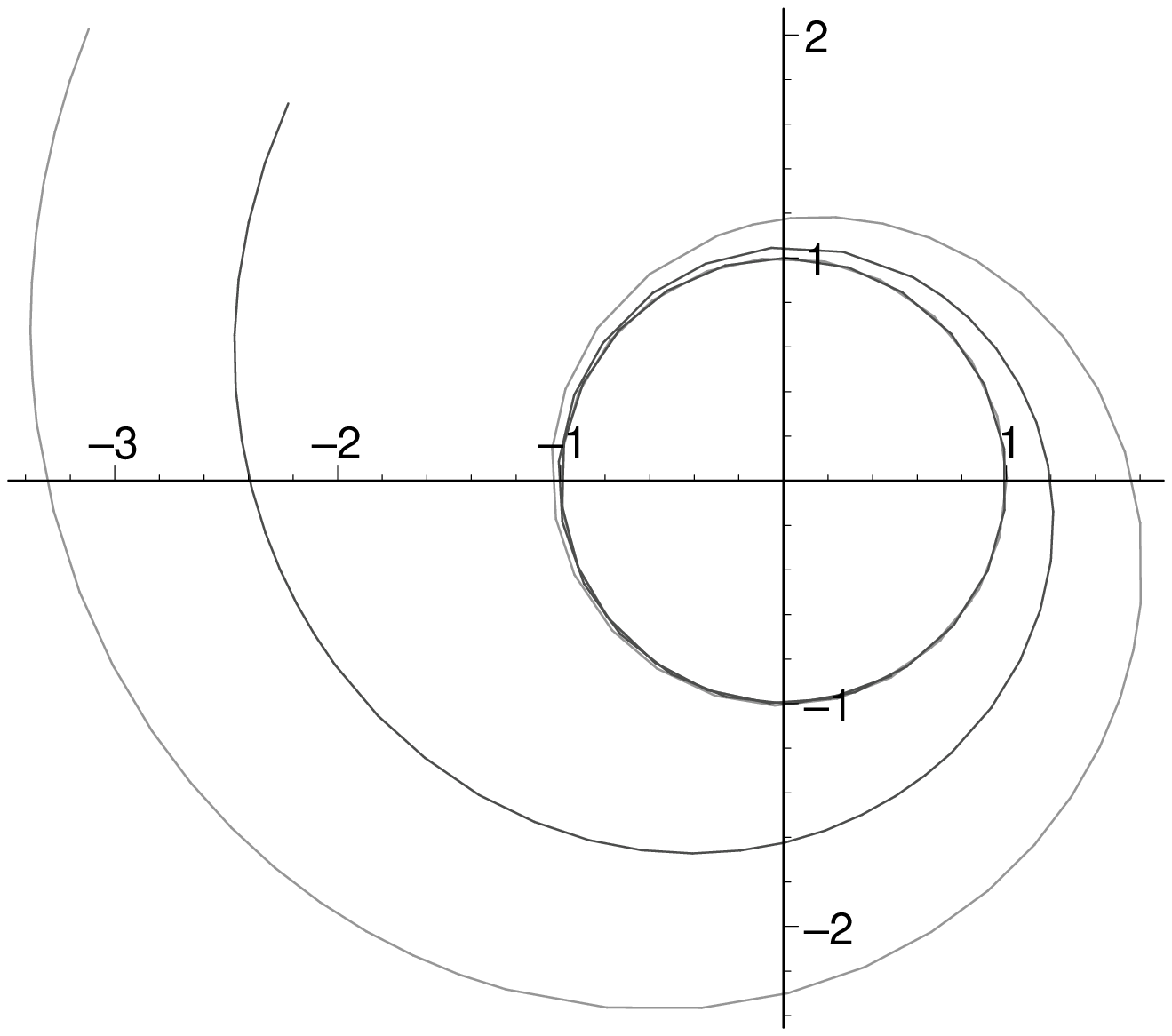}
\caption{}\label{fig1}
\end{center}
\end{minipage}
\begin{minipage}[h]{7cm}
\[
\omega= w_\varrho d\varrho + w_\phi d\phi + w_z dz.
\]
Obviously, $\omega$ vanishes nowhere 
and condition (\ref{notw}) is satisfied. But in 
any neighbourhood of $\{\varrho=1\}$ Euler's multiplier cannot exist. 
Here, the geometry of those manifolds $F$, which are locally the level sets
of a function $f$,
forbid the existence of $f$ and $\lambda$  since they wind 
around the cylinder $\{\varrho=1\}$ (see Fig. \ref{fig1}).
\end{minipage}
\end{figure}
As we can deduce from this example, the geometry of these manifolds, which, 
in general, cannot be represented as submanifolds of $M$ and which we will call
{\it leaves}, has a decisive impact on the existence of $\lambda$ and $f$. 

The geometric term useful in this context is the codim-1-{\it foliation} 
induced by an integrable one-form (see
$\S$2). 
Obviously, in case of global existence of $\lambda$ and $f$, all these leaves are closed subsets
of $M$. Moreover, it turns out that the holonomy group (see $\S$2) has to be trivial.

The first question is, whether those two conditions are sufficient for the global existence
of Euler's multiplier. It turns out that this is only right if we additionally assume the leaves to be compact.

In general, to find a sufficient criterion for global existence, the differentiable structure of $M/{\cal F}$ 
is decisive. In $M/{\cal F}$, points of $M$ lying on the same leaf are identified, in other words,
$M/{\cal F}$ is the space of leaves.
In general, the topology in $M/{\cal F}$ induced by the canonical
projection $\pi: M\to M/{\cal F}$ is non-Hausdorff.

Note that for non-Hausdorff one-dimensional manifolds homeomorphic classes are different form diffeomorphic classes. So, if we only use
topological properties of $M/{\cal F}$ we cannot expect the full regularity of Euler's multiplier.
In order to get a $C^r$ multiplier we have to assume an additional property of ${\cal F}$ involving the differentiable structure. We call this property
{\em regularly $C^r$}. The following example will shed some light on it:
\begin{figure}[h!]
\begin{minipage}[h]{7cm}
\begin{center}
\includegraphics[width=36ex]{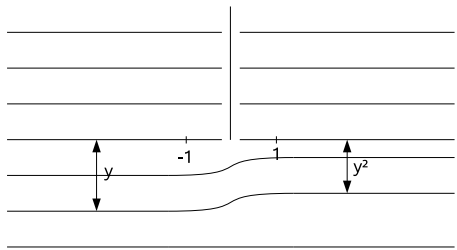}
\caption{An irregular foliation on the slit plane}
\end{center}
\end{minipage}
\begin{minipage}[h]{7cm}
By construction, any $C^1$ function $f$ constant on the leaves would satisfy $f(x,y)=f(-x,-y^2)$, if $x<-1$ and $y<0$. Hence $\frac{\partial f}{\partial y}(x,0)=0$, if
$x<-1$. So $df(x,0)=0$. This prevents the existence of Euler's multiplier although the topological structure is admissible (compare $\S$ 5).  
\end{minipage}
\end{figure}

On the other hand, for Hausdorff one-dimensional manifolds there is no difference between homeomorphic classes and diffeomorphic classes. So, if we
can prove $M/{\cal F}\cong{\bb R}$ in the topological sense, this will be enough to have Euler's multiplier.

In order to reduce $M/{\cal F}$ to its relevant topological properties,
the terminology and the theory of graphs turn out to be  appropriate
means. More precisely, the structure of $M/{\cal F}$ will be displayed
by a configuration of graphs. 
Here, the so called {inseparable} or non-Hausdorff leaves, 
i.e.\ leaves that cannot be separated in the topology of $M/{\cal F}$ play
an important part. They  determine essentially the structure of the graphs.
The construction of the corresponding graphs is possible under 
some weak assumptions on the foliation ${\cal F}$
(cf.\ Def.\ \ref{General assumptions}).
This leads us in a natural way to the concept of the graphical 
configuration induced
by a codim-1-foliation (see $\S$4), which  finally yields an 
equivalent criterion for the
global existence of Euler's multiplier.

In particular, we are able to give a 'classification' of obstructions for the 
global existence of Euler's multiplier, all but one can be expressed in terms
of the topology of $M/{\cal F}$: (1.) non-regularity of ${\cal F}$, (2.) there is at least one non-closed point in $M/{\cal F}$,
(3.) there is too much ramification in $M/{\cal F}$, (4.) there are other
obstructions coming from the local topology of $M/{\cal F}$ or
(5.) obstructions coming from the global topology of $M/{\cal F}$.

For an example for (3.) on $M={\bb R}^2$ we refer to \cite[Ch.\ III, Notes]{ccaln} (found by Wazewski) or \cite{h}.

The authors would like to express their gratitude to Wolf von Wahl for initiating this project.

\section{Foliations and one-forms}

\subsection{Foliations}
Let $M$ be a real $C^r$ manifold of dimension $N$. A maximal $C^r$ atlas 
${\cal F}=(U_\iota,\varphi_\iota)_{\iota\in {\cal J}}$
of $M$ is called a {\it codim-$k$-foliation} if and only if
\begin{enumerate}
\item  For each $\iota\in {\cal J}$ there exist open disks 
$X_\iota\in\real^{N-k}, Y_\iota\in \real^k$ such that 
\[ 
\varphi_\iota(U_\iota)=X_\iota\times Y_\iota,
\]
\item if $(U_{\iota_1},\varphi_{\iota_1}), (U_{\iota_2},\varphi_{\iota_2})
\in {\cal F}$ with $\ U_{\iota_1}\cap U_{\iota_2}\not=\emptyset\ $ then the mapping
\[
\varphi_{\iota_2}\circ\varphi_{\iota_1}^{-1}:\ 
\varphi_{\iota_1}(U_{\iota_1}\cap U_{\iota_2})\to \varphi_{\iota_2}(U_{\iota_1}\cap U_{\iota_2})
\] 
is of the form
\[
\bigl(\varphi_{\iota_2}\circ\varphi_{\iota_1}^{-1}\bigr)(x,y)\ =\ \bigl( \phi_{\iota_1 \iota_2}(x,y), \psi_{\iota_1 \iota_2}(y)\bigr).
\]
\end{enumerate}
The charts $(U_\iota,\varphi_\iota)$ are called {\it foliation charts}.
The leaves $F$ of ${\cal F}$ are defined as follows: Fix any point 
$p\in U_{\iota_0}\subset M$. For any point $q\in M$ we will say $q\in F_p$ if
and only if there exists a finite sequence 
$(\iota_\mu)_{\mu=0,...,m}\subset {\cal J}$ and 
$c_{\iota_0}\in Y_{\iota_0},...,c_{\iota_m}\in Y_{\iota_m}$
such that 
\begin{enumerate}
\item $U_{\iota_{\mu-1}}\cap U_{\iota_{\mu}}\not= \emptyset$, $\ \mu=1,...,m$,
\item $\varphi_{\iota_{\mu-1}}^{-1} \bigl( X_{\iota_{\mu-1}}\times \{ c_{\iota_{\mu-1}}\}\bigr) \cap \varphi_{\iota_{\mu}}^{-1} \bigl( X_{\iota_{\mu}}\times \{ c_{\iota_{\mu}}\}\bigr)\not= \emptyset$, $\ \mu=1,...,m$,
\item $p\in \varphi_{\iota_{0}}^{-1} \bigl( X_{\iota_{0}}\times \{ c_{\iota_{0}}\}\bigr)$, $\ q\in \varphi_{\iota_{m}}^{-1} \bigl( X_{\iota_{m}}\times \{ c_{\iota_{m}}\}\bigr)$ 
\end{enumerate}
It is clear, that by 
$
p\overset{\cal F}{\sim}q\, :\Longleftrightarrow q\in F_p
$
an equivalence relation is defined. 
Any equivalence class $F=F_p=[p]_{\overset{\cal F}{\sim}}$, 
$p$$\in$$M$, will be called a leaf of ${\cal F}$. \\
The space of leaves is denoted by 
$M/{\cal F}:=M/{\underset{\mbox{\large $\sim$}}{\mbox{\tiny ${\cal F}$}}}$,
and 
$\pi: M\to M/{\cal F}, p\mapsto F_p$ is the canonical projection. 
In this article, we consider $M/{\cal F}$ as a topological space 
endowed with the topology
induced by $\pi$. The projection $\pi$ then becomes a 
continuous and open mapping.
\\
In what follows, for any leaf $F$ of ${\cal F}$ we will write
$F\in M/{\cal F}$ if we consider $F$ as a point in $M/{\cal F}$ and,
deviating from the exact definition of the foliation,
$F\in {\cal F}$ if $F$ is considered as a subset of $M$.\\
At this occasion, we point out that, of course, the sets $\varphi_{\iota}^{-1} \bigl( X_{\iota}\times \{ c\}\bigr)$, $c\in Y_{\iota}$, are submanifolds of
dimension $N$$-$$k$, but in general the leaves $F\in\cal F$ are not necessarily
submanifolds; they are immersed manifolds. 
For this and  other theorems in the theory of foliations used in this 
article we refer to \cite{ccaln}.\\
For any subset $U\subset M$ we define the {\it saturation} of $U$ with respect
to ${\cal F}$ by ${\rm sat}\ (U):=\pi^{-1} \bigl(\pi (U)\bigr)$.
Any subset $U\subset M$ is called {\it saturated} (with respect to ${\cal F}$)
if and only if ${\rm sat}\ (U)=U$, i.e.,
for each $p\in U$ holds $F_p\subset U$.

\

We call the foliation {\em transversely orientable}, if the normal bundle given by $N_{{\cal F},x}:=T_xM/T_xF_x$ is orientable. The transversely orientable $\codim$-1-foliations
are exactly those given by integrable one-forms. 
%
%
\subsection{The holonomy group}
We now give a short sketch of the construction of holonomy. Fix a leaf
$F$ and a point $p\in F$. Now consider a closed path $\gamma:[0;2\pi]\seq F$
with $\gamma(0)=\gamma(2\pi)=p$. Now, if we choose a local transversal $T$
through $p$ and cover $\gamma([0;2\pi])$ by a finite number of 
foliation charts $(U_\iota,\varphi_\iota)_{\iota=1,...,k}$, 
we can choose points $p_i\in\gamma([0;2\pi])$ and local transversals $T_i$ through
$p_i$ such that $T_i\cup T_{i+1}\subset U_\iota$ for a $\iota\in\{1,...,k\}$.
By going into the foliation charts we are now able to choose a path
$\gamma_t:[0;2\pi]\seq F_t$ for a $t\in T$ sufficiently near to $p$, which 
intersects each $T_i$. 

Indeed, it can be shown, that $\gamma_t(2\pi)$ is independent of the choices of
$p_i$, $T_i$ and the chosen paths $\gamma_t$. Moreover, $\gamma_t(2\pi)$
only depends on the homotopy class of $\gamma$. Hence we get a map
$$h:\pi_1(F,p)\seq {\rm Aut}(T,p),$$
where $\pi_1(F,p)$ is the fundamental group of $F$ with respect to $p$ and 
${\rm Aut}(T,p)$ denotes the germs of automorphisms of $T$ which fix $p$. 
We call
$$\hol(F,p):=h\bigl(\pi_1(F,p)\bigr)$$
the holonomy of the leaf $F$. Indeed, $\hol(F,p)$ is independent of the
base point $p$, so that we will denote this by $\hol(F)$.
For more details of the construction we refer to \cite{ccaln}. 

%

%
%
\subsection{The foliation of an integrable one-form}
Let $M$ be a real manifold of dimension $n$$+$$1$, $\omega$ a $C^r$ one-form on $M$.
Thus, for each $p\in M$ 
\[
\omega(p):\quad T_{p}M\to \real
\]
is a linear mapping.
Setting 
\[
P(p):={\rm ker}\ \omega(p),
\]
where $\omega(p)\not= 0$, an $n$-plane field $P$ is defined on $M$.\\
Frobenius' Theorem yields that $P$ is induced by a codim-1-foliation on $M$ if
and only if $P$ is completely integrable, which is equivalent to the condition
\[
\omega \wedge d\omega =0
\]
on $M$ (see \cite[Ch.II $\S$ 4, Appendix $\S$ 3 Th.2]{ccaln}).\\
Thus, each completely integrable $C^r$ one-form $\omega$ on $M$ with $\omega(p)\not=0$ in $M$ induces
a $C^r$ codim-1-foliation ${\cal F}={\cal F}(\omega)$ on $M$, where for $p\in F$,
$F\in {\cal F}$ 
\[
T_{p} F\, =\, {\rm ker}\ \omega(p).
\]
Now, considering the corresponding foliation charts of $M$, the {\it local} existence of 
Euler's multiplier is clear.

Moreover, the global existence of Euler's multiplier easily implies that $\hol(F)=\{1\}$ for all leaves $F$.

Also note that the foliation given by an integrable, nowhere vanishing one-form $\omega$ is transversely orientable: By using a partition of unity one
can easily construct a $C^r$ vector field $X$ satisfying $\omega(X)\equiv 1$.

Finally, we want to mention the Godbillon-Vey class of ${\cal F}(\omega)$ (see \cite{gv}), if $r\ge\dim M\ge 3$. 
Since $\omega\wedge d\omega=0$, we can find a one-form $\eta$ such that
$d\omega=\omega\wedge\eta$. Now the Godbillon-Vey class is defined
$$gv(\omega):=[-\eta\wedge d\eta]\in H^3(M,{\bb R})$$
as a de Rham class. Taking the cohomological class instead of the form itself has the effect, that $gv(\omega)=gv(f\omega)$ for any non-vanishing $f\in C^{r}(M)$.
So it depends only on the foliation, not on the one-form. If $\omega$ allows for Euler's multiplier, $gv(\omega)=0$.

\section{Transverse fibrations}

%

%

%

For the sake of simplicity, in this section all manifolds and maps are considered to be $C^\infty$ unless stated otherwise.
A fibre bundle $\phi:M\seq B$ with $M,B$ being manifolds is called {\em transverse
to the foliation ${\cal F}$ on $M$}, if and only if
$$\phi|F:F\seq B$$
is an unramified, surjective covering map for every leaf $F$ and therefore, in particular,
$$T_pM=T_pF\oplus T_pM_{\phi(p)},$$
where $p\in F$ and $M_{\phi(p)}:=\phi^{-1}(\phi(p))$ is the fibre of $\phi$ over $\phi(p)$. In this situation we will also say, that
${\cal F}$ {\em is transverse to} $\phi$. 
In our context, we assume, that $\dim M=n+1$ and $\dim B=n$ and hence, that ${\cal F}$ is a codimension one foliation.

For such $\phi$ we can define a {\em global} version of holonomy in the following way.
Given a closed path $\alpha:[0;2\pi]\seq B$ on $B$ and a point 
$y\in\phi^{-1}(\alpha(0)):=T$, we can 
look at the lifting $\tilde\alpha:[0;2\pi]\seq F_y$ on the leaf $F_y$ containing
$y$ with $\tilde\alpha(0)=y$. Mapping $y\mapsto\tilde\alpha(2\pi)$ induces
a $C^r$ diffeomorphism $\phi_\alpha:T\seq T$. In fact, $\phi_\alpha$ only
depends on the homotopy class of $\alpha$. So we get a representation
$${\bs h}:\pi_1(B)\seq \diff^r(T),$$
where $\pi_1(B)$ is the fundamental group of $B$ and $\diff^r(T)$ the group of $C^r$ automorphisms of $T$.
We call the image of ${\bs h}$ the global holonomy group (with respect to $\phi$)
of ${\cal F}$ and denote this by $\hol({\cal F}/B)$. For details see again \cite{ccaln}.

As $\phi|F:F\seq B$ is a covering map, we can embed $\pi_1(F)\subset\pi_1(B)$ 
and look at the restricted map
$${\bs h}|F:\pi_1(F)\seq \diff^r(T),$$
whose image we call $\widetilde{\hol}(F)$, on the one hand and at the
ordinary holonomy map
$$h_F:\pi_1(F)\seq \hol(F)$$
on the other hand, both of them are group homomorphisms. 
By definition of ${\bs h}$ follows that 
$$\ker({\bs h}|F)\subset \ker(h_F)$$
and hence there is a surjective homomorphism
$$\widetilde{\hol}(F)=\pi_1(F)/\ker({\bs h}|F)\seq\pi_1(F)/\ker(h_F)=\hol(F)\.$$
We now proved that $\hol(F)$ is a quotient of a subgroup of $\hol({\cal F}/B)$
for every leaf $F$.

So $\hol({\cal F}/B)=1$ implies $\hol(F)=1$ for all leaves $F$, but the converse
is not true.

\begin{Theorem}\label{transversefibrations1} Let $\phi:M\seq B$ a fibre bundle and 
${\cal F}$ a transverse codim-1-foliation, both of class $C^r$. Then are equivalent:
\begin{enumerate}
{\rm \item \hspace{0ex}}There is a $C^r$ function $f:M\seq{\real}$ such that $df(x)\not=0$ in $M$ and ${\cal F}={\cal F}(df)$.
{\rm \item \hspace{0ex}}$\hol({\cal F}/B)=1$ and the fibres of $\phi$ are diffeomorphic to ${\real}$.
\end{enumerate}
In particular, $M\cong B\times{\real}$ if one (and hence both) of the conditions
is satisfied.
\end{Theorem}

\begin{proof}
'$\Rightarrow$': Let $F$ be a leaf and $T$ a fibre of $\phi$, meeting in
a point $p\in M$. Then we obtain a non-vanishing
vector field $X$ along $T$ satisfying
$$df|_T(X)>0\.$$
Now let $\gamma:[0,y]\seq T$ be a differentiable path in $T$ such that
$$d\gamma(t)(1)=X(\gamma(t)),$$
By Stokes' Theorem $f(\gamma(y))-f(\gamma(0))>0$.
So $T$ is diffeomorphic to ${\real}$ and intersects each leaf exactly once.
This implies $\hol({\cal F}/B)=1$, in particular.

'$\Leftarrow$': Let $T$ be a fibre of $\phi$. The condition $\hol({\cal F}/B)=1$
implies that $T$ intersects each leaf $F$ in exactly one point.
This means that there is a differentiable map $\psi:M\seq T$.
Now take a
diffeomorphism $\delta:T\seq{\real}$ and define $f:=\delta\circ\psi$.
Then the leaves of ${\cal F}$ are the level sets of $f$. The construction easily yields $df(x)\not=0$ everywhere.
%
\end{proof}

Let us recall that the relative tangent bundle $T_{M|B}$ is defined by
$$T_{M|B,x}=di_x(T_xM_{\phi(x)}),$$
where $M_{\phi(x)}$ is again the fibre over $\phi(x)$, $i:M_{\phi(x)}\seq M$ the natural embedding and $di_x$ the derivative of $i$ at $x$,
embedding the tangent bundle of the fibre into the tangent bundle of $M$. The $i$-th Betti number of a manifold $B$ is denoted by
$b_i(B)=\dim H^i(B,{\bb R})$. 

\begin{defn}{}Let $M$ be a manifold and $\phi:M\seq B$ a fibre bundle. We call $M$ {\em $\phi$-orientable} if and only if $T_{M|B}$ is orientable.\end{defn}

\begin{Theorem}Let $M, B$ be manifolds, $\phi:M\seq B$ a fibre bundle such that $\dim M=\dim B+1$ and
$M$ is $\phi$-orientable. If $b_2(B)=0$,
then there is a transverse and transversely orientable foliation ${\cal F}$ and a Riemannian metric on $M$ such that $T_{M|B}=T_{\cal F}^\perp$
and all fibres of $\pi$ are geodesics.

If  $B$ is simply connected, then every transverse
foliation ${\cal F}$ is transversely orientable.
\end{Theorem}

\begin{proof}By general theory \cite[Ch.\ V, \S 4]{ccaln} any two transverse foliations ${\cal F},{\cal F}'$ with the same global holonomy are equivalent, i.e.\, 
there is some $H\in \diff(M)$ such that $\phi\circ H=\phi$ and takes leaves of ${\cal F}$ to leaves of ${\cal F}'$. If ${\cal F}'={\cal F}(\omega)$ for 
an integrable one-form $\omega$, then $H^*\omega$ is an integrable one-form and ${\cal F}={\cal F}(H^*\omega)$. Furthermore $\hol({\cal F}/B)=1$ for all
transverse foliations ${\cal F}$, if $\pi_1(B)=1$. So it is sufficient to construct an integrable one-form $\omega$
such that ${\cal F}(\omega)$ is transverse to $\phi$. 

By assumption, $M$ is $\phi$-orientable and hence there is a non-vanishing vector field $X\subset T_{M|B}$. For any $p\in M$ we can find a chart $V\subset B$ such that
$\phi^{-1}(V)\cong C\times V$, $C$ being the fibre. After a choice of a section $S$ of $\Phi$ over $V$ we can introduce the fibre coordinate $t$ of $x\in\phi^{-1}(V)$ by
$x=\alpha_{X}(t)$; here $\alpha_\cdot$ denotes the flow of the subscript vector field starting in $S(\phi(x))$. If $C\cong{\bb R}$, then $dt$ is
a well defined one-form on $\phi^{-1}(V)$. If $C\cong S^1$, then we can achieve
$$\alpha_X(t)=\alpha_X(t+1)$$ 
by multiplying $X$ with a nonvanishing global function, independently of the choice of $S$. 
Given this property, the one-form $dt$ is well defined on $\phi^{-1}(V)$.

If we now think of $M=\bigcup_i C\times V_i$ to be covered by such charts, in the intersection $C\times(V_i\cap V_j)$ we have
$$dt^{(i)}=dt^{(j)}-\phi^*\eta_{ij}$$
on $\phi^{-1}(V_i\cap V_j)$
for a closed one-form $\eta_{ij}\in{\cal A}^1_B(V_i\cap V_j)$ depending on the sections $S_i$ and $S_j$. 
Let us denote by ${\cal K}\subset{\cal A}^1_B$ the sheaf of closed one-forms on $B$. If we knew $H^1(B,{\cal K})=0$, then we could find closed one-forms
$\eta_i\in{\cal A}^1_B(V_i)$ such that $\eta_{ij}=\eta_i-\eta_j$ on $V_i\cap V_j$. Thus
$$\omega:=dt^{(i)}+\phi^*\eta_i$$
would define a global closed one-form with $\omega(X)=1$, yielding a transverse and transversely orientable foliation.

In order to compute $H^1(B,{\cal K})$ note that
$$0\seq{\cal K}\seq{\cal A}^1_B\stackrel{d}{\seq}{\cal A}^2_B\stackrel{d}{\seq}\dots$$
is an acyclic resolution abbreviating the standard acyclic resolution of ${\bb R}$ by one step. So we find
$H^1(B,{\cal K})=H^2(B,{\bb R})=0$.

After an appropriate choice of a Riemannian metric $g_B$ on $B$ we can construct a Riemannian metric
$$g:=\omega\otimes\omega+\phi^*g_B$$
on $M$. The geodesic equation for the fibre is $\Gamma^{i}_{11}=0$, first for all $i>1$, but then also for $i=1$, since geodesics are parametrised proportional to arc length. This translates to $\omega_{i,1}=0$, a condition satisfied since $\omega$ is closed and $\omega_1=1$.
Finally, if $\tau\in T_{{\cal F},x}$ we compute
$$g(\frac{\partial}{\partial t},\tau)=\omega(\frac{\partial}{\partial t})\omega(\tau)+\phi^*g_B(\frac{\partial}{\partial t},\tau)=0+0=0,$$
so $T_{M|B}=T_{\cal F}^\perp$.

%
\end{proof}

\begin{ex}{}Every circle bundle (with structure group $U(1)\subset Aut(S^1)$, that is)
$\pi:M\seq B$ is $\pi$-oriented.
So, if $B$ is a manifold such that $H^2(B,{\bb Z})$ is a torsion group, every non-zero element of $H^2(B,{\bb Z})$ induces a non-trivial circle bundle over $B$ satisfying the
assumptions of the theorem. For example, real projective space of dimension $\ge 2$ or every complex Inoue surface is such a base manifold $B$. 
\end{ex}

Up to now, in all cases, for which we proved ${\cal F}={\cal F}(df)$,
we had $M/{\cal F}\cong{\real}$.

\section{Regular foliations induced by integrable one-forms}

%

%


We want to formulate regularity conditions in terms of global generatedness. For this purpose recall a relative version of this. 

\begin{defn}For any sheaf
${\cal E}$ on a topological space $M$ and continuous map $\psi:M\seq N$ to a topological space $N$ we say that
${\cal E}$ is {\em $\psi$-generated}, if for all $x\in M$, open sets $U\subset M$ with $x\in U$ and $f\in{\cal E}(U)$ there exist $V\subset U, W\subset N$ open 
and $\tilde f\in{\cal E}(\psi^{-1}(W))$ such that $x\in V\cap\psi^{-1}(W)$ and $\tilde f|V=f|V$. 
We say ${\cal E}$ is {\em globally generated}, if $\psi$ can be chosen to be the constant map.
\end{defn}

For foliations there is natural sheaf to consider.

\begin{defn}Let $M$ be a $C^r$ manifold and ${\cal F}$ be a $C^r$ foliation. For any open set $U\subset M$ we denote
$$C^r_{\cal F}(U):=\{f\in C^r(U)| f \mbox{ is constant on the leaves of }{\cal F}|U\}\.$$
\end{defn}

\begin{defn}\label{big} \begin{lenumerate}
\item \llabel{sim} For leaves $F,F'$ of a $C^r$ foliation ${\cal F}$ on $M$ we say $F$ is {\em infinitely close} to $F'$, if for any saturated open sets
$U\supset F, U\supset F'$ holds $U\cap U'\not=\emptyset$. The smallest equivalence relation generated by this property
is denoted by $\sim$, its quotient by $G$ and $p:M/{\cal F}\seq G$ the quotient map. Finally, let $S:=\{F\in{\cal F}|\,\,p^{-1}(p(F))\not=\{F\}\}$ denote
the {\em non-Hausdorff leaves} of ${\cal F}$.
\item The foliation ${\cal F}$ is called
{\it of finite type}, if
\begin{lenumerate}
\item \llabel{eins} each leaf $F\in{\cal F}$ is closed,
\item  \llabel{zwei}$p(S)$
is a discrete subset of $G$, 
and for each $g\in G$ the set $p^{-1} (g)$ is finite.
\end{lenumerate}
\item For a $C^r$-foliation ${\cal F}$ on $M$ we call $\tilde M_{\cal F}:=(N_{\cal F}\setminus(\{0\}\times M))/{\bb R}^+$ the {\em ${\cal F}$-oriented double cover of $M$};
here ${\bb R}^+$ acts by multiplication on the fibres. We denote the natural projection by $\alpha:\tilde M_{\cal F}\seq M$ and obtain a foliation $\tilde{\cal F}:=\alpha^*{\cal F}$ with projection maps $\tilde\pi:\tilde M_{\cal F}\seq\tilde M_{\cal F}/\tilde{\cal F}$ and $\tilde p:\tilde M_{\cal F}/\tilde{\cal F}\seq\tilde G$. 
The foliation $\tilde{\cal F}$ is transversely orientable.
\item A $C^r$-foliation ${\cal F}$ on $M$ is called {\em regularly $C^r$},
if the sheaves $\tilde\pi_*C^s_{\tilde{\cal F}}$ are globally generated for all $1\le s\le r$. If $r=\infty$, we call ${\cal F}$ regularly $C^\infty$, if
$\tilde\pi_*C^\infty_{\tilde{\cal F}}$ is globally generated.

\item A $C^r$-foliation ${\cal F}$ on $M$ is called {\em almost regularly $C^r$}, if the sheaves $\tilde\pi_*C^s_{\tilde{\cal F}}$ are $\tilde p$-generated for all $1\le s\le r$;
if $r=\infty$, we impose only $\tilde\pi_*C^\infty_{\tilde{\cal F}}$ to be $\tilde p$-generated.
\end{lenumerate}
\end{defn}

The seemingly complicated definition of regularity stems from the effort to separate the combinatorial data of ${\cal F}$ from its analytic. If we would have replaced
the ${\cal F}$-oriented double cover of $M$ simply by $M$, non-regularity could have occurred just for combinatorial reasons, e.g. the foliation in Figure \ref{fig2}(c)
would be non-regular, independently of the analytic behaviour of the leaves near $F_1', F_2', F_3'$. The ${\cal F}$-oriented double cover in this case would have six
equivalent non-Hausdorff leaves (with micrograph being a hexagon, cf. next chapter) and the well-meaning reader may agree that the extension property given in the definition
of regularity holds in this case. 

We will see later that non-regularity is a crucial obstruction for a foliation of finite type to be characterized by combinatorial data. So 
we want to inquire into conditions for a foliation of finite type to be regular.


Our main question is: When is a foliation induced by a one-form regular? For the sake of simplicity we restrict ourselves to the case of $C^\infty$ regularity.
The corresponding results for $C^r$-regularity are similar and obtained in the same way. 

\begin{Lemma}\label{reg1}Let $\omega$ be an integrable $C^\infty$ one-form without zeroes on the $C^\infty$ manifold $M$. If the foliation ${\cal F}(\omega)$ is
of finite type and 
there exists $\lambda\in  C^\infty(M)$ such that $\omega\wedge d\lambda=d\omega$, then ${\cal F}(\omega)$ is regularly $C^\infty$.  
\end{Lemma}

\begin{proof}Let $X$ be a $C^\infty$ vector field on $M$ such that $\omega(X)\equiv 1$. This can easily be constructed using a partition of unity. Now we consider
the flow $\alpha_{t}$ of $h  X$ for a $h \in C^\infty(M)$. We want to choose $h $ in such a way that leaves will be mapped to leaves (where
defined), i.e. for any leaf $F$ and its inclusion $i_F:F\inj M$ we want $i_F^*\alpha_t^*\omega=0$ for all admissible $t$. Furthermore no leaf shall be fixed, i.e.
$h$ is wanted without zeroes.
Differentiating this equation by $t$
and looking at $t=0$ yields
$${\cal L}_{h  X}\omega\wedge\omega=0,$$
where ${\cal L}$ denotes the Lie derivative. 

On the other hand, ${\cal L}_{h  X}\omega\wedge\omega=0$ implies ${\cal L}_{h  X}\omega=\kappa\omega$ for some $\kappa\in C^\infty(M)$ and we compute
$$\frac{d}{dt}i_F^*\alpha_t^*\omega=i_F^*\frac{d}{dt}\alpha_t^*\omega=i_F^*\alpha_t^*{\cal L}_{h  X}\omega=\kappa i_F^*\alpha_t^*\omega,$$
hence
$$i_F^*\alpha_t^*\omega=C\exp(\kappa t)i_F^*\omega=0\.$$
So our derived equation is indeed equivalent to the leaf preserving property of the flow.

By elementary rules (see e.g. \cite{lang}) we compute
\begin{eqnarray*}{\cal L}_{h  X}\omega&=&h {\cal L}_X\omega+\omega(X)dh \\
&=&hd(\omega(X))+hC_Xd\omega+dh \\
&=&hC_Xd\omega+dh .\end{eqnarray*}
Here $C_.$ denotes the contraction map.
Note furthermore that $C_Xd\omega\wedge\omega=C_X(d\omega\wedge\omega)-\omega(X)d\omega=-d\omega$. So our leaf conserving equation reads
$$dh \wedge\omega=hd\omega,$$
or, assuming $h>0$ and defining $\lambda:=\log h$,
$$d\lambda\wedge\omega=d\omega\.$$

Now we can choose foliation charts such that $h  X$ is the tangent vector associated to the non-leaf coordinate function 
if and only if $h $ is nowhere zero. These
foliation charts patch together to give an almost regular $C^\infty$ foliation due to the non-vanishing of $h $: 
Let $f\in C^\infty_{\cal F}(U)$ for a saturated 
open set $U$, $F\subset U$ a leaf and $F'$ infinitely close to $F$. We find 
$\varepsilon>0, x\in F, x'\in F'$ such that $$\psi:(-\varepsilon,\varepsilon)\seq M, t\mapsto\alpha_t(x)\mbox{ and }\psi':(-\varepsilon,\varepsilon)\seq M, t\mapsto\alpha_t(x')$$ are well-defined.
Due to the facts that $F$ and $F'$ are infinitely near and $h  X$ is leaf-preserving, 
we have that the points $\psi(t)$ and $\psi'(t)$ are on the same leaf whenever
$\psi'(t)$ is on a leaf contained in $U\cap\psi((-\varepsilon,\varepsilon))$. 
Let $U'$ be the saturation of $\psi'((-\varepsilon,\varepsilon))$. Since $d\psi'(0)=h (x')X(x')\not=0$, we may assume that $U'$ is open.
Furthermore, since ${\cal F}$ is transversely orientable and all leaves are closed
we may assume that $\pi\circ\psi$ and $\pi\circ\psi'$ are injective with non-vanishing differentials. So the function
$$\tilde f:U\cup U'\seq{\bb R}, z\mapsto\left\{\begin{array}{cl}f(z)&\mbox{, if }z\in U\\f(\psi(t))&\mbox{, if }z=\psi'(t)\in U'\end{array}\right.$$
is well-defined, an element of $C^\infty_{\cal F}(U\cup U')$ and extends $f$ beyond $F'$. By the way of construction, these extensions glue together to yield
an extension $\tilde f'$ of $f$  on a saturated open set $\tilde U\supset [F]_\sim\cup U$, since ${\cal F}$ is of finite type.
The same assumption about ${\cal F}$ allows us to extend $\tilde f'|\tilde V$ to a function $\tilde f\in C^\infty_{\cal F}(M)$ vanishing outside $\tilde U$ for a shrunk
saturated open set $U\cup[F]_\sim\subset\tilde V\subset\tilde U$. 
So ${\cal F}$ is regularly $C^\infty$.
%
%
\end{proof}

This implies immediately that ${\cal F}(\omega)$ is regular, if it is of finite type and $\omega$ is closed.
The reader will agree that the proof of Lemma \ref{reg1} is also valid to prove almost regularity, if the condition on ${\cal F}$ to be of finite type is relaxed
appropriately.
The first order differential equation characterizes that some multiple of $\omega$ is closed:

\begin{Lemma}\label{diff1}$\omega\wedge d\lambda=d\omega$ has a global solution $\lambda\in C^\infty(M)\iff \exists\,h\in C^\infty(M)$ positive
such that $h\omega$ is closed.\end{Lemma}

\begin{proof}For any positive function $h\in C^\infty(M)$ and global solution $\lambda$
of $\omega\wedge d\lambda=d\omega$, the function $\tilde\lambda:=\lambda-\log h$ solves $\tilde\omega\wedge d\tilde\lambda=d\tilde\omega$ for $\tilde\omega:=h\omega$. 
Hence we choose $h:=\exp(\lambda)$ in order to achieve $\tilde\lambda=0$. If $h\omega$ is closed, then $0=hd\omega+dh\wedge\omega$, hence $\log h$ is
a global solution of the equation. 
\end{proof}

The $C^\infty$ integrable one-form $\omega$ induces a canonical element $[\omega]\in H^1(M,C^\infty_{\cal F})$ via the following construction. Let $\omega=f_idg_i$ on $U_i$ as constructed by Frobenius' Theorem and $\lambda_i:=\log|f_i|$.
This is a solution to the differential equation in Lemma \ref{diff1}, so $d(\lambda_i-\lambda_j)\wedge\omega=0$ on
$U_i\cap U_j$, i.e. $\phi_{ij}:=\lambda_i-\lambda_j\in C^\infty_{\cal F}(U_{ij})$. If the covering
$U_i$ is chosen fine enough, the functions $\phi_{ij}$
define an element $[\omega]\in H^1(M,C^\infty_{\cal F})$ independent of the choices made.
 
\begin{Proposition}Let $M$ be a $C^\infty$ manifold, $\omega$ a nowhere vanishing integrable $C^\infty$ one-form and ${\cal F}={\cal F}(\omega)$ of finite type. 
For the properties
\begin{enumerate}
\item $[\omega]=0\in H^1(M,C^\infty_{\cal F})$,
\item ${\cal F}$ is induced by a {\em closed} one-form,
\item ${\cal F}$ is regularly $C^\infty$
\end{enumerate}
holds: $(a)\imp (b)\imp (c)$.

\end{Proposition}

\begin{proof}

'$(a)\imp (b)$':
We construct $\phi_{ij}$ like above.
Property $(a)$ now yields $\phi_i\in C^\infty_{\cal F}(U_i)$ such that $\lambda_i+\phi_i=\lambda_j+\phi_j$ on 
$U_i\cap U_j$. So we obtain a global solution $\lambda\in C^\infty(M)$ for $d\lambda\wedge\omega=d\omega$.

The implication $(b)\imp (c)$ is already proved.
\end{proof}  

If ${\cal F}(\omega)$ is given by a closed one-form, then $gv(\omega)=0$, of course.
Note that the property $gv(\omega)=0$ is much weaker than the property that ${\cal F}$ is given by a closed one-form. 
In case of $M={\bb R}^3$ a smooth version of the first example of the introduction (leaves winding around a cylinder)
cannot be given by a closed one-form, but $gv(\omega)=0$, since $b_3(R^3)=0$. 

Regularity will turn out to be the main assumption in order to enable us to decide the existence of Euler's multiplier from graphical data constructed in the next section (cf. proofs of Lemma \ref{semiglobal} and Theorem \ref{Hauptsatz}). 

\section{The graphical configuration of a foliation and the existence of Euler's multiplier}

%

%

%

\begin{figure}[ht]
\begin{center}
\vspace{-10ex}
\includegraphics[width=93.6ex]{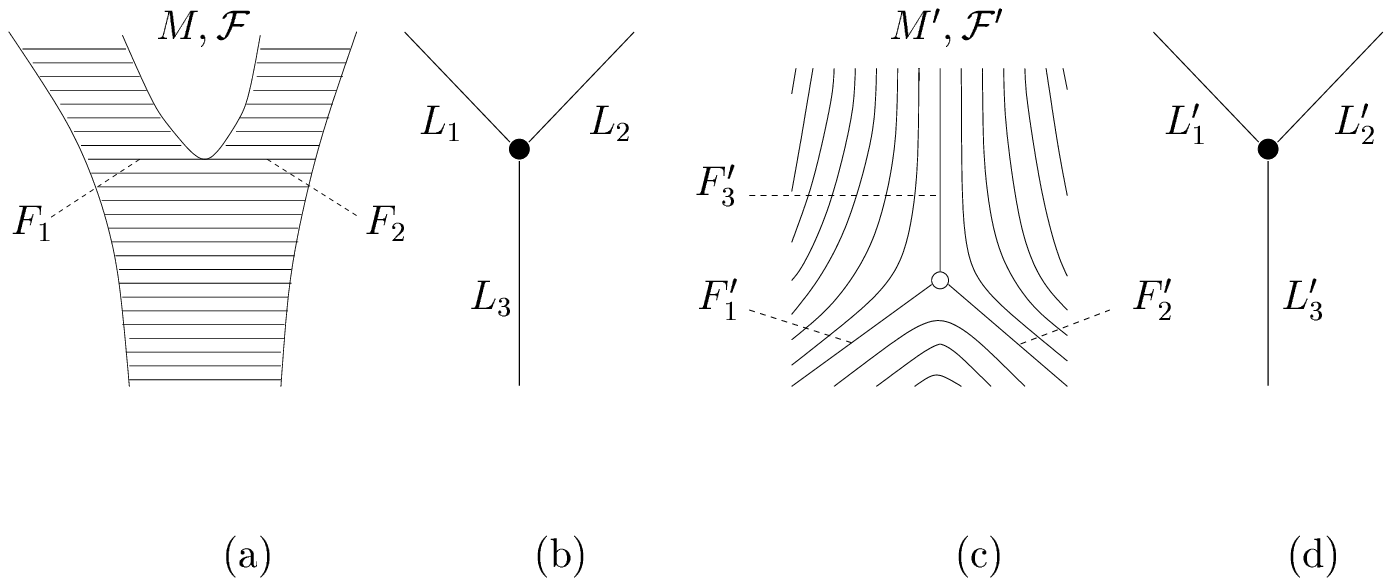}
\vspace{-99ex}
\caption{Examples of foliations with non-Hausdorff leaves\label{fig2}}
\end{center}
\end{figure}

\begin{figure}[ht]
\vspace{-13ex}
$\mbox{\hspace{9ex}}$
\includegraphics[width=80ex]{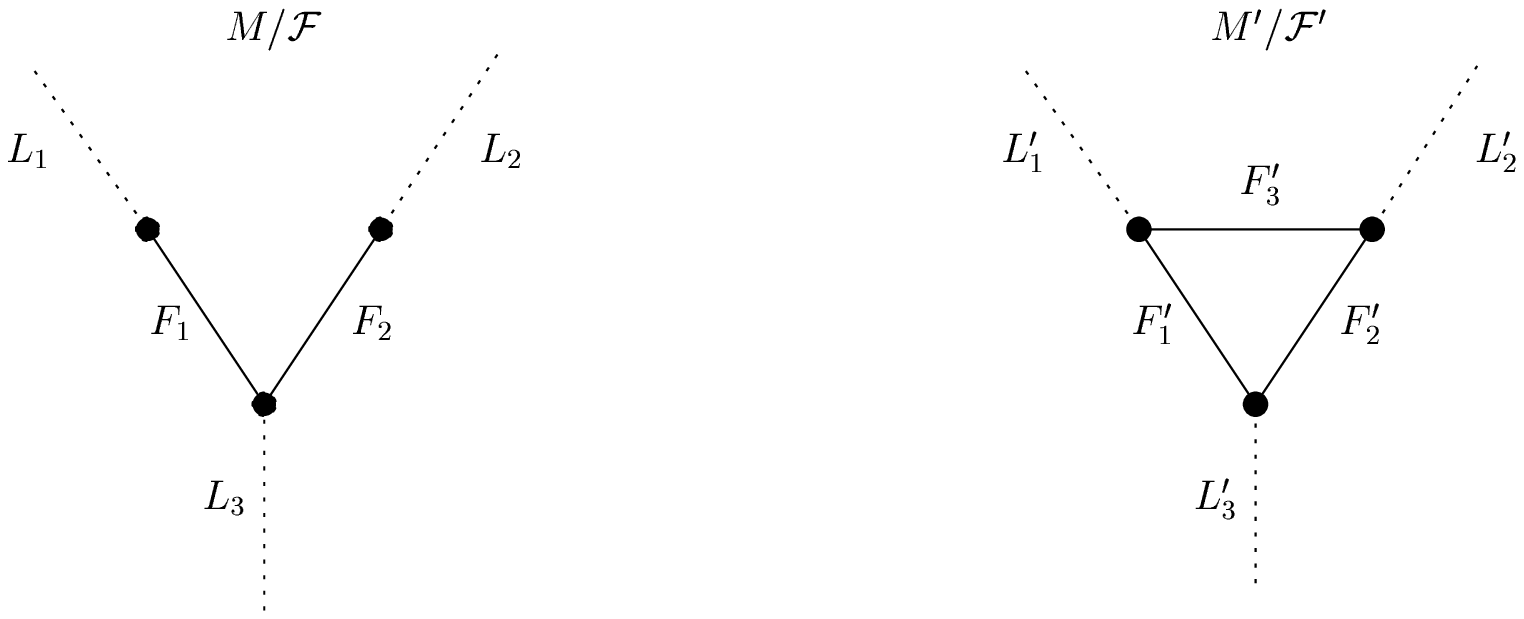}
\vspace{-78ex}
\caption{Corresponding micrographs to the examples above\label{fig3}}
\end{figure}

If $M/{\cal F}$ is a one-dimensional real manifold diffeomorphic to $\real$,
it is now clear that there is some function $f\in C^1(M,\real)$ with
\begin{equation}\label{f}
\forall x\in M:\quad df(x)\not= 0,\quad \mbox{ and } \quad {\cal F}={\cal F}(df).
\end{equation}
Roughly speaking,
in what follows, we are going to decompose $M/{\cal F}$ into manifolds
diffeomorphic to $\real$
and investigate consistency at the connecting points. The global existence
of $f$ as in (\ref{f}) then depends on a certain geometric and topological
configuration of $M/{\cal F}$ which can be described by graphs.\\
The components of $M/{\cal F}$ which are manifolds diffeomorphic
to $\real$ are
identified with the edges of a graph, in the following called the {\it macrograph}
of $M/{\cal F}$. The vertices of this macrograph correspond to the irregular
points of $M/{\cal F}$, i.e., points that are non-Hausdorff-points or endpoints of $M/{\cal F}$.
For example, the leaves $ F_1$  and $ F_2$ in 
the foliation ${\cal F}$
in Fig.\ \ref{fig2} (a) or the leaves $ F'_1, F'_2, F'_3$ of ${\cal F}'$ in
Fig.\ \ref{fig2} (c) are such non-Hausdorff-points. (The marked point in Figure \ref{fig2}(c)
where the leaves $F_1',F_2',F_3'$ meet together does not belong to $M$.)
The corresponding (local) macrographs of these ``bifurcations'' of ${\cal F}$ or
${\cal F}'$ are the same and
will look like
in Fig.\ \ref{fig2} (b) or Fig.\ \ref{fig2} (d), respectively. The macrographic configuration of a foliated manifold is not alone
decisive for the existence of $f$ 
with the properties mentioned above.
The foliation ${\cal F}'$
in Fig.\ 2 has the same macrographic configuration as ${\cal F}$, but obviously,
${\cal F}$ admits such a function $f$ as in (\ref{f}), 
whereas ${\cal F}'$ does not so. 
\\
Therefore, the behaviour of the foliation ${\cal F}$ in the irregular points
of $M/{\cal F}$ is to be taken into account. It will be described by so-called
{\it micrographs}.
\\
The growth of $f$, if it exists, will later be denoted by an orientation of the
edges of the macrograph.
The micrograph in a vertex of the macrograph regulates the ``traffic'', i.e.,
tells us along which ways in $M/{\cal F}$ the function $f$ is increasing.
Therefore the macroedges beginning or ending in some macrovertex are identified
with the vertices in the micrograph. The microedges are the irregular points of
$M/{\cal F}$ and connect the two macroedges from which they are accessible.\\
For example, the foliations in Figure \ref{fig2} have in the macrovertices,
displayed as bullets in Figure \ref{fig2}(b) and (d), the
micrographs as in Figure \ref{fig3}.
The system consisting of the macrograph $\Gamma$ and the micrographs 
$\gamma_\sigma$ in the vertices $\sigma\in\Gamma$, integrated in $\Gamma$,
is called the graphical configuration of $M/{\cal F}$ and is decisive for
the global existence of Euler's multiplier.\\
It will turn out that two criteria, one on the micrographs and one on the
macrograph, are necessary and, together with some topological conditions
on the leaves of ${\cal F}$, also sufficient for the global existence 
of Euler's multiplier.\\
In a first step, supposing only the first criterion, we will prove the
{\it semiglobal} existence, i.e., if all micrographs admit the globalization
of $f$, then for each $p\in M$ there is a  neighbourhood 
$\tilde U\subset M$ of $F_p$ and, should the occasion arise, all other leaves
belonging to the same macrovertex as $F_p$ (see Def.\ \ref{big}\lref{sim}, Def.\ \ref{Makrograph}), saturated 
with respect to ${\cal F}$, such that in
$\tilde{U}$ Euler's multiplier can be defined.
In this case, if the orientation of one macroedge is chosen, 
the micrographs of
the macrovertices  connected to this macroedge determine the 
orientation of the neighbouring macroedges. \\
If this condition of semiglobal existence is fulfilled, then the 
orientability of the macrograph with respect to the regulation
of the micrographs and the topological properties of the 
oriented macrograph decide about the global existence.\\
In particular, if $M$ is simply connected, 
the macrograph and the micrographs of ${\cal F}$ are also
simply connected. By the latter property semiglobal existence is guaranteed, but this even implies global existence by the simply connectedness of the
macrograph. 
\\
\begin{figure}[h!]
\begin{center}
\includegraphics[height=25ex]{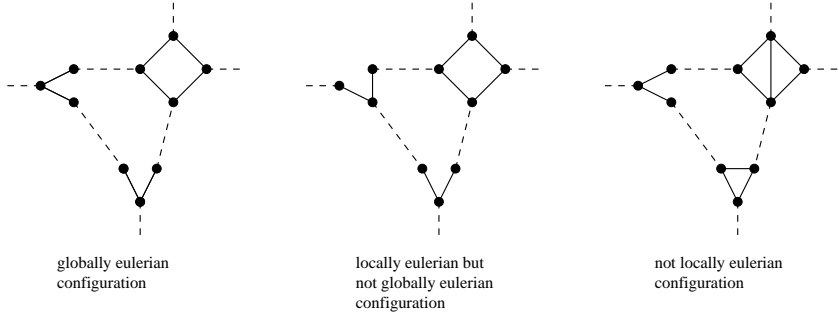}
\caption{Examples of locally or globally eulerian configurations\label{configur1}}
\end{center}
\end{figure}
In Figure \ref{configur1} the graphs of three foliations corresponding
to the same macrographs are shown. The edges of the macrograph are
represented
by dotted lines, the edges of the micrographs of the
macrovertices by full lines. 
The terms {\it locally eulerian} and {\it globally eulerian} are defined
in \ref{eulerian}.
\\
The ideas mentioned above shall be stated more precisely in the following.
%
%
%
\subsection{Construction of the graph of a foliation}
%
%
%
Having discussed foliations of a certain type in the previous chapter, we want to go into greater detail here and
therefore give them a name.
\begin{defn} \label{General assumptions} 
\rm
Let $M$ be a real $C^r$-manifold of dimension $n$$+$$1$, ${\cal F}$
a codim-1-foliation of class $C^r$ on $M$. The foliation ${\cal F}$ is called
{\it graphical} if and only if it is of finite type and regularly $C^r$. 
\end{defn}
If a Jordan Brouwer separation theorem holds on $M$, we can show that the two conditions in the definition of finite type
foliations are not independent. 
This is done in $\S 6$.\\
Note that the finiteness of $p^{-1}(g)$ for each $g\in G$ implies that
$G$ is Hausdorff.\\[0.5ex]
In what follows, we assume ${\cal F}$ to be a graphical foliation on a
real manifold $M$ according to Definition \ref{General assumptions}. 
%
%
%
%
%
%
%
%
%
\begin{figure}[h!]
\begin{minipage}[h]{7cm}
\begin{defn} \label{endpoints}
\rm (Endpoints)
Let $e\in M/{\cal F}$. We say that $e$ is an {\it endpoint} if and only if there
exists some $a\in M/{\cal F}$ such that
for each continuous injective mapping $w: [0,1]\to M/{\cal F}$
with
$w(0)=a$ and $e\in w([0,1])$ holds $w(1)=e$.
Moreover, let $E$ denote the set of endpoints in $M/{\cal F}$.
\end{defn}
For example, the leaf $F$ in Figure \ref{fig4} is  such  an endpoint.
\end{minipage}
\begin{minipage}[h]{7cm}
\begin{center}
\includegraphics[width=20ex]{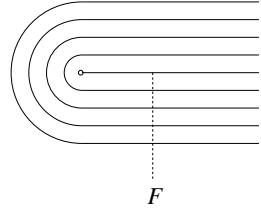}
\end{center}
\caption{Example of an endpoint\label{fig4}}
\end{minipage}
\end{figure}
Now we are in a position to construct the macrograph of the
foliation ${\cal F}$.
%
%
%
%
%
%
%
%
%
%
\begin{construction} \label{Makrograph}  
\rm (Macrograph) 
(a) {\it Vertices}: If $(S/\sim)\cup (E/\sim)=\emptyset$ take
an arbitrary $\sigma\in G$ and set 
${\cal V}\Gamma^0:=\{\sigma\}$,
else let
\[
{\cal V}\Gamma^0:=(S/\sim)\cup (E/\sim).
\]
Let us define recursively:
(i) If there is a
continuous injective mapping $w:[0,1)\to G$ 
with $w(0)=\sigma_1\in {\cal V}\Gamma^{j}$ 
and $w((0,1))\subset G
\setminus {\cal V}\Gamma^{j}$ and 
the limit $\displaystyle\lim_{t\to 1}w(t)$ does not exist, then define
\[
{\cal V}\Gamma^{j+1}:={\cal V}\Gamma^{j}\cup
\bigcup_{k\in{\nat}}w\Bigl(1-\frac 1k\Bigr).
\]
(ii) If there is an injective mapping $w:[0,1)\to G$ 
with $w(0)=\sigma_1\in {\cal V}\Gamma^{j}$ 
and $w((0,1))\subset G
\setminus {\cal V}\Gamma^{j}$ and $\displaystyle\lim_{t\to 1} w(t)=w(0)=\sigma_1$ then
define
\[
{\cal V}\Gamma^{j+1}:={\cal V}\Gamma^{j}\cup
w\Bigl(\frac12\Bigr).
\]
As there can only be a countable number of ranges 
$w((0,1))\subset G$, by choice of a counting it can be obtained that
\[
{\cal V}\Gamma:=\bigcup_{j\in{\nat}} {\cal V}\Gamma^{j}
\]
has the property that $w(1):=\lim_{t\to 1}w(t)$ does exist for 
every continuous, injective map 
$w:[0,1)\to G$ with $w(0)=\sigma_1\in {\cal V}\Gamma$ and $w((0,1))
\subset G\setminus {\cal V}\Gamma$, where $w(1)\not=w(0)$.\\
Now ${\cal V}\Gamma$ is called the set of {\it vertices}
of the macrograph $\Gamma$. Note that $S\subset p^{-1}({\cal V}\Gamma)$, 
but, in general, 
there are also noncritical vertices (cf.\ Remark \ref{artificialvertices}).
\\[0.5ex]
(b) {\it Edges}: To define the edges of the macrograph 
first we set
\[
\begin{split}
{\bf W}:=\Big\{
w:[0,1]\to M/{\cal F}\ \Big\vert\ w &\mbox{ is continuous, } w\big\vert_{[0,1)}
\mbox{ is injective, }\\
&w(0), w(1)\in p^{-1}({\cal V}\Gamma), 
w((0,1))\cap p^{-1}({\cal V}\Gamma)=\emptyset
\Big\}.
\end{split}
\]
Now we say $w\sim^{*} w'$ for $w,w'\in {\bf W}$ if and only if
\[
p\Bigl( w\bigl([0,1]\bigr)\Bigr)\, =\, p\Bigl( w'\bigl([0,1]\bigr)\Bigr)
\]
For each $L\in {\bf W}/\sim^{*}$ set an edge $L$ connecting
the corresponding points $p(w(0))$ and $p(w(1))$ in 
${\cal V}\Gamma$,
i.e., 
let ${\cal E}\Gamma:={\bf W}/\sim^{*} $ be the set of edges of $\Gamma$.\\
Now we can define the macrograph $\Gamma=({\cal V}\Gamma, 
{\cal E}\Gamma)$ with the corresponding mapping
\begin{equation}\label{kappamacro}
\kappa:\ {\cal E}\Gamma \to {\it S}^{2} {\cal V}\Gamma,\
L\mapsto \{ p(w(0)), p(w(1))\},\quad\, \mbox{ where }\ [w]_{\sim^{*}}=L.
\end{equation}
%
\end{construction}
%
%
%
%
%
For $x\in M$ and $L\in {\cal E}\Gamma$ we will write
\begin{equation}\label{xinL}
x\in L\, :\Leftrightarrow\, \exists w\in{\bf W}: [w]_{\sim^*}=L,\ \pi(x)\in w((0,1)).
\end{equation}
%
%
%
%
%
%
%
%
\begin{rem}\label{artificialvertices} \rm
The construction in \ref{Makrograph} (a) guarantees that each edge of 
$\Gamma$ connects two different vertices. This is due to the formal 
definition of graphs in literature. The addition (i) 
is necessary to avoid edges running to infinity and is illustrated in 
Figure \ref{fig5}. 
\begin{figure}[h!]
\begin{minipage}[h]{7.0cm}
\begin{center}
\includegraphics[height=13ex]{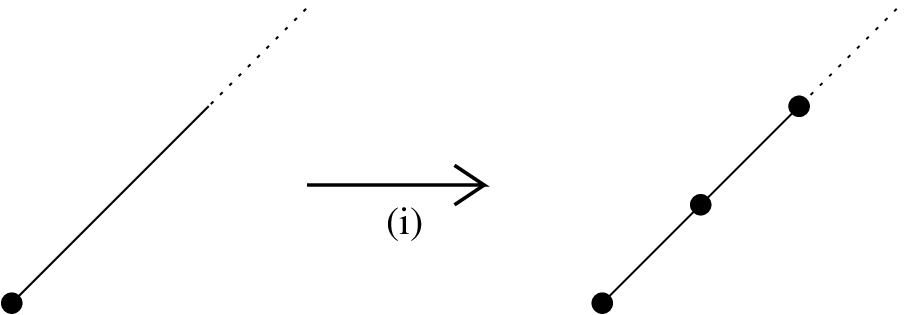}
\caption{construction (i) in \ref{Makrograph}\label{fig5}}
\end{center}
\end{minipage}
\begin{minipage}[h]{7.0cm}
\begin{center}
\includegraphics[height=13ex]{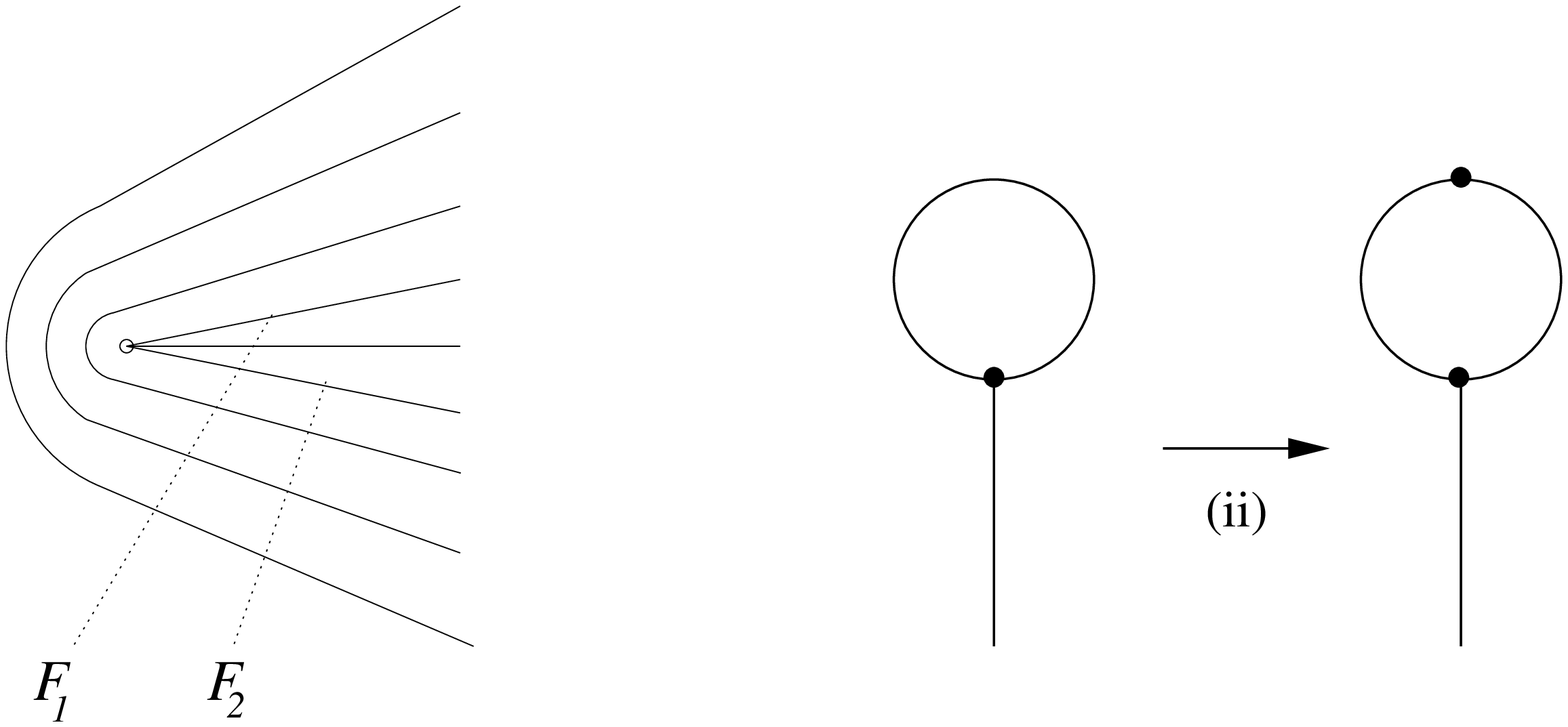}
\caption{construction (ii) in \ref{Makrograph}\label{fig6}} 
\end{center}
\end{minipage}
\end{figure}
The addition (ii) suppresses loops in $\Gamma$ by adding a 'noncritical'
vertex. An example of a foliation generating such a loop and the 
corresponding macrograph is shown in Figure \ref{fig6}.
%
%
%
Here, the leaves
$F_1$ and $F_2$ are non-Hausdorff and generate the first vertex. 
Note that these 'artificial' vertices have trivial inner structure and are
therefore noncritical, cf.\ Remark \ref{viersieben}. 
\end{rem}
Now we investigate the macrovertices and construct the 
corresponding micrographs.
%
%
%
%
%
%
%
%
%
\begin{construction} \label{Mikrograph} 
\rm (Micrograph) 
Let $\sigma\in {\cal V}\Gamma \setminus (E/\sim)$,
$p^{-1}(\sigma)=\big\{ s_1,...,s_N\big\}\subset M/{\cal F}$. 
Let 
${\cal E}\Gamma^{\sigma}:=\{ L_1,...,L_\nu\}$ be the set of macroedges $L$
with $\sigma\in \kappa (L)$, where $\kappa$ is defined in (\ref{kappamacro}).
 \\
(a) {\it Vertices}: Let 
\[
{\cal V}\gamma_\sigma:=\big\{ L_1,...,L_\nu\big\}
\]
be the set of vertices of the micrograph $\gamma_\sigma$ of the
macrovertex $\sigma$.\\
(b) {\it Edges}:
We say for $s\in p^{-1}({\cal V}\Gamma)$ and 
$L\in {\cal E}\Gamma$
\[
s\in {\it Lim}\ (L)\ :\Longleftrightarrow \ 
\exists w\in {\bf W}: 
\ [w]_{\sim^{*}}=L, w(1)=s.
\]
Now we will show that for each leaf $s_j\in p^{-1}(\sigma)$ there are
exactly two different macroedges $L$ with $s_j\in {\it Lim} (L)$:\\
Since $s=s_j\in p^{-1}(\sigma)$ is no endpoint, there is
a continuous bijective mapping $w: [0,1]\to M/{\cal F}$ with
$s\in w((0,1))$. Let $\tau\in (0,1)$ with $w(\tau)=s$.
Extending and rescaling 
$w\vert_{[0,\tau]}$ and $w\vert_{[\tau,1]}$
in an appropriate way we can obtain
two different macroedges $[w_1]_{\sim^*}$ and $[w_2]_{\sim^*}$ with
$s\in {\it Lim} ([w_1]_{\sim^*})\cap {\it Lim} ([w_2]_{\sim^*})$.
Thus there are at least two macroedges connected to $s$.
Now let $x\in s\subset M$ be any fixed point and $(U,\varphi)$ 
a foliation chart of $(M,{\cal F})$, where $x\in U$, 
$D\subset\real^n$, $I=(-\epsilon, +\epsilon)$  and
$\varphi=(\varphi_1,\varphi_2): U\to D\times I$ is a bijective $C^r$ mapping.
Now consider the commutative diagram
$$\xymatrix{
M\supset\mbox{\hspace{-8ex}} & 
U\ar[d]^\pi\ar[r]^{\varphi} & D\times I\ar[d]^{pr}\\
M/{\cal F}\supset\mbox{\hspace{-6ex}} & V\ar[r]^{i_V} & I
},$$
where $V=\pi(U)$ and $i_V$ is, by definition of the local manifold
structure of $M/{\cal F}$, a diffeomorphism. 
Without restriction we can assume $\varphi_2 (x)=0$
and therefore we have 
\[
i_V (w_j((0,1))\cap V)= (-\epsilon,0)\, 
\mbox{ or }\, i_V (w_j((0,1))\cap V)= (0,+\epsilon),\quad\, j=1,2,
\]
if $\epsilon >0$ is sufficiently small. But 
$i_V (w_1((0,1))\cap V)=i_V (w_2((0,1))\cap V)$
implies already $[w_1]_{\sim^*} =[w_2]_{\sim^*}$.
Thus, there are exactly two macroedges $L$ satisfying 
$s\in {\it Lim} (L)$.\\
Now we set ${\cal E}\gamma_\sigma:=\{ s_1,...,s_N\}$. For $j=1,...,N$
and 
$s\in {\it Lim} (L_{j_1})\cap {\it Lim} (L_{j_2})$, 
where $L_{j_1}\not=L_{j_2}$ set the microedge $s$ 
connecting the microvertices $L_{j_1}$ and $L_{j_2}$. 
Now the micrograph $\gamma_\sigma$ is given by 
$\gamma_\sigma=({\cal E}\gamma_\sigma, {\cal V}\gamma_\sigma)$
with the corresponding mapping
\begin{equation}\label{kappamicro}
\kappa_\sigma :\quad {\cal E}\gamma_\sigma \to {\it S}^{2} {\cal V}\gamma_\sigma
,\quad s_j\mapsto \{ L_{j_1}, L_{j_2}\}
.
\end{equation}
\end{construction}
%
%
%
%
%
%
\begin{rem}\label{viersieben}  \rm
Examples of micrographs are shown in Figure \ref{fig3}. The micrograph
of a noncritical macrovertex has trivial structure which is 
displayed in Figure \ref{fig7}.
\begin{figure}[ht]
\begin{center}
\includegraphics[width=20ex]{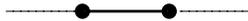}
\caption{A micrograph of a noncritical macrovertex\label{fig7}}
\vspace{-3.5ex}
\end{center}
\end{figure}
\end{rem}
%
%
%
%
%
%
%
%
\begin{defn}\label{graphconf} {\rm (Graphical configuration)} Let ${\cal F}$ be a graphical
foliation on a real manifold $M$ according to Definition 
\ref{General assumptions}, $\Gamma$ the macrograph of $M/{\cal F}$ and
$\gamma_\sigma$ the micrographs of the macrovertices 
$\sigma\in {\cal V}\Gamma$ according to Constructions \ref{Makrograph},
\ref{Mikrograph}.
Then $\delta (F):=\bigl(\Gamma, (\gamma_\sigma)_\sigma\bigr)$ is called
the {\it graphical configuration} of ${\cal F}$.
\end{defn}
%
%
%
%
%
%
%
%

The following is just a refinement of \cite[Ch.\ 4, Lemma 5]{ccaln}.

\begin{Lemma}\label{hol}
Let ${\cal F}$ be a foliation with only closed leaves on the manifold $M$. 
Then for every leaf $F$ of ${\cal F}$ the holonomy is 
$\hol(F)\in\{1,{\Bbb Z}_2\}$. Furthermore, if ${\cal F}$ is graphical and $\hol(F)={\Bbb Z}_2$ then
$F\in E$, i.e. $F$ is an endpoint of
$\delta({\cal F})$, where $\delta({\cal F})=(\Gamma,(\gamma_\sigma)_\sigma)$ is 
the graphical configuration of
the foliation ${\cal F}$.
\end{Lemma}
The following notion of oriented and bipartite graphs will turn out to
be appropriate for the criteria of existence of Euler's multiplier.
%
%
%
%
%
%
%
%
\begin{defn}  \label{orientation}
{\rm (Orientation)} \rm 
Let ${\cal G}$ be the macrograph $\Gamma$ (or any other graph).
Then a mapping 
\[
\uparrow\ :\, \, {\cal E}{\cal G} \to {\cal V}{\cal G}\times {\cal V}{\cal G}
\]
is called orientation of ${\cal G}$
if and only if $\,\ {\rm f}\ \circ \uparrow\ =\kappa,\ $ where
\[
{\rm f}:\, {\cal V}{\cal G}\times {\cal V}{\cal G}\to S^2{{\cal V}{\cal G}},\,
(\sigma_1,\sigma_2)\mapsto\{ \sigma_1,\sigma_2\}.
\]
\end{defn}
%
%
%
%
%
%
%
%
\begin{defn}\label{bipartite}  
\rm (bipartite Graph) 
Let ${\cal G}=\bigl({\cal VG},{\cal EG}\bigr)$ be a micrograph $\gamma_\sigma$ or any other graph with corresponding mapping $\kappa$ 
(cf.\ (\ref{kappamicro})). The graph
${\cal G}$ is called {\it bipartite} if and only if there exist two subsets
${\cal V}_+{\cal G}$, ${\cal V}_-{\cal G}\subset {\cal VG}$ with the following properties
\begin{enumerate}
\item ${\cal VG}={\cal V}_+{\cal G}\cup {\cal V}_-{\cal G}$,\quad 
${\cal V}_+{\cal G}\cap {\cal V}_-{\cal G} =\emptyset$,
\item $\forall s\in {\cal EG}\ \ \exists\ L^+\in {\cal V}_+{\cal G}, \
L^-\in {\cal V}_-{\cal G}:
\ \ \kappa(s)=\{L^+, L^-\}$.
\end{enumerate}
The second condition means that only vertices of different classes are
connected by edges of ${\cal G}$.
\end{defn}
The micrographs of the first and the second example in Fig.\ \ref{configur1}
are bipartite, whereas two of the micrographs in the right example
of Fig.\ \ref{configur1} are not bipartite.

%

%

%
\subsection{The foliation of a graphical configuration}

We want to revert the construction of the graphical configuration. For this purpose we generalize Definition \ref{graphconf}:

\begin{defn}A {\em graphical configuration} is a tuple $\left(\Gamma,(\gamma_\sigma)_\sigma\right)$ with a graph $\Gamma$, graphs $\gamma_\sigma$ for every
$\sigma\in{\cal V}\Gamma$ such that ${\cal E}\Gamma^\sigma={\cal V}\gamma_\sigma$.
\end{defn}

For the definition of ${\cal E}\Gamma^\sigma$ see Construction \ref{Mikrograph}.

Now we can prove

\begin{Theorem}\label{constr}
Let $\delta=\left(\Gamma,(\gamma_\sigma)_\sigma\right)$ be a graphical configuration without endpoints. Then there is closed one-form $\omega\not=0$ on a
surface $M$ (maybe with zeroes), such that the zero locus of $\omega$
is smooth and has no zero-dimensional components,
the leaves of the so defined foliation
${\cal F}\left(\omega\right)$ are diffeomorphic to ${\real}$ and
$\delta=\delta\left({\cal F}\left(\omega\right)\right)$.
\end{Theorem}

\begin{proof}
Let $\left(\Gamma,\uparrow\right)$ be an orientation of $\Gamma$, that means
$${\cal E}\Gamma\stackrel{\uparrow}{\seq}{\cal V}\Gamma\times{\cal V}\Gamma\stackrel{{\rm f}}{\seq}S^2{\cal V}\Gamma,$$ where
f is the forgetful map (forgetting the order) and ${\rm f}\circ\uparrow=\kappa$, cf. Definition \ref{orientation}.
Let us denote an element $m\in{\cal V}\Gamma\times{\cal V}\Gamma$ by $\left(m_1,m_2\right)$, and construct as follows. Let for $\sigma\in{\cal V}\Gamma,
L\in{\cal E}\Gamma$ denote
$$Lim^\sigma\left(L\right):=\left\{s\in{\cal E}\gamma_\sigma| L\in\kappa_{\sigma}\left(s\right)\right\},$$
where $\kappa_{\sigma}$ denotes the boundary map of $\gamma_\sigma$
(cf.\ (\ref{kappamicro})).
For every $L\in{\cal E}\Gamma$ define
$$I_{L}:=[0,1]\times\left(0,1\right)\setminus\left(\left\{0\right\}\times\left\{\left.\tfrac ik\right| i=1,...,k-1\right\}\cup\left\{1\right\}\times
\left\{\left.\tfrac il\right|\ i=1,...,l-1\right\}\right),$$
if $Lim^{\uparrow_1\left(L\right)}L=\left\{s_{11}\left(L\right),...,s_{1k}\left(L\right)\right\}$ and $Lim^{\uparrow_2\left(L\right)}L=
\left\{s_{21}\left(L\right),...,s_{2l}\left(L\right)\right\}$. This means, the $s_{ab}(L)$ are the limit points of $L$, where $a=1,2$ stands for the two boundary
vertices of $L$, ordered according to the chosen orientation.

Now take
$$M=\bigcup_{L}I_{L}/\sim',$$
where $\sim'$ is defined as follows: If $s\in Lim^{\uparrow_i\left(L\right)}L\cap
Lim^{\uparrow_j\left(L'\right)}L'$
for one (and then a unique) choice $i,j\in\left\{1,2\right\}$ and $s=s_{ia}\left(L\right)=s_{jb}\left(L'\right)$ then
$$ I_{L}\ni\left(i-1,\frac{a-x}{k}\right) 
\sim'\left(j-1,\frac{b-x}{l}\right)\in I_{L'}$$
for $x\in\left(0,1\right)$. This means, that the intervals in $I_L$ resp $I_L'$ corresponding to $s$ are identified.

The so defined topological space $M$ can easily be made to a differentiable manifold, because after the choice of $\sigma$,
the edges $\left\{L,L'\right\}=\kappa_{\sigma}\left(s\right)$ are unique and therefore a small neighbourhood of every point of $M$ is topologically an open
set in ${\real}^2$. The differentiability of the transition maps is just a formal calculation.
The foliation ${\cal F}$ is
given by $x=const$ if $\left(x,y\right)\in I_{L}$. Since the coordinates
of the $I_L$ can differ at most in the sign, there is a global one-form
$\omega$, which is closed and can locally be written
$$\omega=\sin(2\pi x)\ dx\.$$
Hence, $\omega$ defines the foliation ${\cal F}$
even in the zero locus of $\omega$.
\end{proof}

%
%
%
%
%
\begin{figure}[h!]
\begin{minipage}[h]{6.0cm}
\begin{center}
\includegraphics[width=30ex]{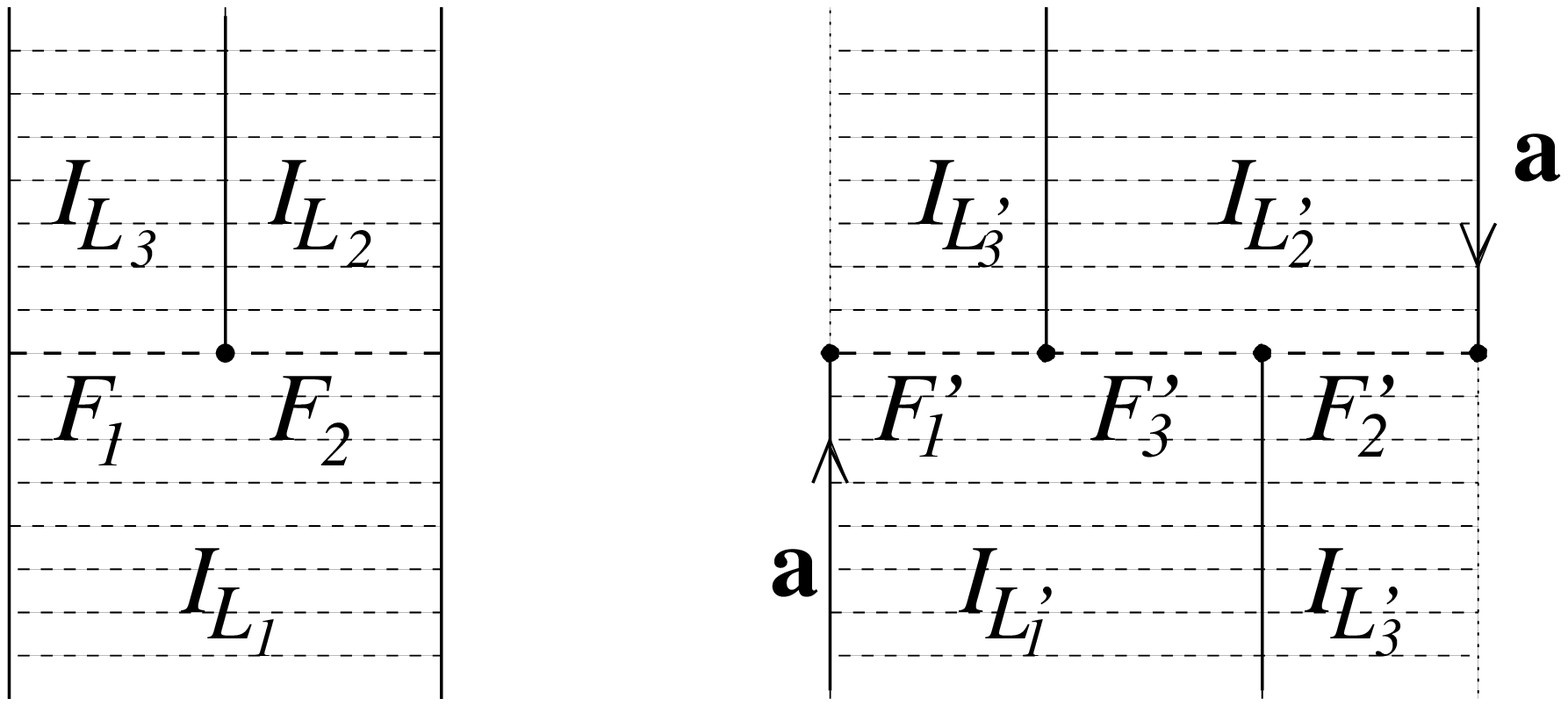}\quad
\caption{
\label{fig8}} 
\end{center}
\end{minipage}
\begin{minipage}[h]{8.0cm}
In Fig.\ \ref{fig8} we illustrated the rather technical construction of
Theorem \ref{constr} for the graphical configurations of Fig.\ \ref{fig3} resp.\
Fig.\ \ref{fig2}. In the second case, the right and left boundaries ``${\rm a}$'' are to
be identified in converse direction. Hence the resulting manifold equals
a M\"obius surface without three certain rays. By a look at the manifold
in Fig.\ \ref{fig2}(c) we see in particular, that
the correspondence of graphical configurations and the pairs $(M,{\cal F})$
is not reversible.
\end{minipage}
\end{figure}

%

%

%
\subsection{The existence of Euler's multiplier}
Now we arrive at the main result of this paragraph which is an
equivalent criterion on the graphical configuration $\delta({\cal F})$
for the local or the global existence, respectively,
of Euler's multiplier.
%
%
%
%
%
%
%
%
\begin{defn} \label{eulerian}
\rm
(Eulerian configuration)
Let ${\cal F}$ be a graphical $C^r$ codim-1-foliation on $M$,
 $\delta({\cal F})=\bigl(\Gamma, (\gamma_\sigma)_\sigma\bigr)$
the graphical configuration of ${\cal F}$ with respect to the construction
in \ref{Makrograph}, \ref{Mikrograph}.
\begin{enumerate}
\item  The configuration $\delta({\cal F})$ is called
{\it locally eulerian} if and only if $E=\emptyset$ and
for each $\sigma\in {\cal V}\Gamma$ the corresponding micrograph
$\gamma_\sigma$ is bipartite (cf.\ Definitions \ref{endpoints},
\ref{bipartite}).
\item The configuration $\delta({\cal F})$
is called {\it globally eulerian} if and only if the following three
conditions hold:
\begin{enumerate}
\item\ilabel{globallyeulerian} $\delta({\cal F})$ is locally eulerian.
\item The bipartitions 
$
{\cal V}\gamma_\sigma\, =\, {\cal V}_+\gamma_\sigma \cup {\cal V}_-\gamma_\sigma
$
of $\sigma\in {\cal V}\Gamma$
can be chosen such that, using the notation of \ref{Mikrograph}
and \ref{bipartite}, for $\sigma,\sigma'\in {\cal V}\Gamma$
\[
{\cal E}\Gamma^\sigma \cap {\cal E}\Gamma^{\sigma'}\, =\, \bigl(
 {\cal V}_+\gamma_\sigma \cap {\cal V}_-\gamma_{\sigma'}\bigr)
\cup \bigl(
 {\cal V}_-\gamma_\sigma \cap {\cal V}_+\gamma_{\sigma'}\bigr).
\]
\item There are no cycles
$
\sigma_1 L_1 \sigma_2 L_2 \sigma_3 ... \sigma_N L_N \sigma_1
$
with $L_i\in {\cal V}_+\gamma_{\sigma_i}$ for all $i=1,...,N$ or
$L_i\in {\cal V}_-\gamma_{\sigma_i}$ for all $i=1,...,N$, where $N\in\nat$.
\end{enumerate}
\end{enumerate}
\end{defn}
At this occasion, we want to point out that Definition \ref{eulerian}
has nothing to do with the usual definition of an {\it eulerian graph}
in literature.
\begin{rem}\label{bem}
\rm Obviously, the definition of a global
eulerian configuration in Definition
\ref{eulerian} \iref{globallyeulerian} is equivalent to the following: \\
The configuration $\delta({\cal F})$ is
{\it globally eulerian} if and only if there is an orientation
\[
\uparrow\ : {\cal E}\Gamma \to {\cal V}\Gamma\times
{\cal V}\Gamma, \ L\mapsto \bigl( \uparrow_1 (L), \uparrow_2 (L)\bigr)
\]
of $\Gamma$ (cf.\ Def.\ \ref{orientation}) such that
the following three
conditions hold:
\begin{enumerate}
\item[]
\begin{enumerate}
\item $E=\emptyset$.
\item For each $\sigma\in {\cal V}\Gamma$ the micrograph
$\gamma_\sigma$ is bipartite such that, using the notation of Definition
\ref{bipartite},
\[
{\cal V}_+{\gamma_\sigma} = \big\{ L\in {\cal V}\gamma_\sigma \big\vert
\uparrow_1 (L) =\sigma\big\},\quad
{\cal V}_-{\gamma_\sigma} = \big\{ L\in {\cal V}\gamma_\sigma \big\vert
\uparrow_2 (L) =\sigma\big\}.
\]
\item \iilabel{drei} The orientation $\uparrow$ induces a well-defined partial
ordering ``$\le$'' on ${\cal V}\Gamma$ by
\[
\begin{split}
\sigma <\sigma'\ :\Longleftrightarrow\ &\exists N\in\nat,\
\tau_1,...,\tau_N
\in {\cal V}\Gamma, \ L_1,...,L_{N-1}\in {\cal E}\Gamma:\\
& \tau_1=\sigma,\ \tau_N=\sigma',\ 
\uparrow (L_j)=(\tau_{j},\tau_{j+1}),
\mbox{ \footnotesize $j=1,...,N$$-$$1$},
\end{split}
\]
where $\sigma <\sigma'$ means $\sigma \le \sigma'$ but $\sigma \not= \sigma'$.
\end{enumerate}
\end{enumerate}
\end{rem}
%
Examples of locally eulerian or globally eulerian configurations
are illustrated in Fig.\ \ref{configur1}.\\
%
%
%
In order to formulate our criterion for semiglobal existence of
Euler's multiplier, we need some technical preparations:
\begin{defn}\label{Further notations} \rm
A curve $\alpha: [a,b]\to M$ of class $C^r$ is called {\it transversal},
if the mapping $\pi\circ\alpha: [a,b] \to M/{\cal F}$ is injective and
\[
\forall t\in [a,b]: \quad \frac{d}{dt}\bigl( \pi\circ\alpha\bigr) (t)\not = 0.
\]
\end{defn}
\begin{rem}\label{tildeU}  \rm
Let $ U\subset G = (M/{\cal F})/\sim$ an open set
and $p:M/{\cal F}\to G$ the canonical projection
(cf.\ Def.\ \ref{big}\lref{sim}).
Then
$
\tilde U:=\pi^{-1}\bigl(p^{-1}(U)\bigr)
$
is a saturated set in $M$ satisfying
\[
F\in {\cal F},\
F\cap \tilde U\not=\emptyset \,\, \Rightarrow \, \bigcup_{F'\in [F]_{\sim}}
\mbox{\hspace{-2ex}} F'\ \subset\ \tilde U,
\]
i.e.\ if $F\cap \tilde U\not= \emptyset$
for any leaf $F$
then $\tilde U$ is a common
saturated neighbourhood of all leaves $F'\in {\cal F}$
with $F'\in [F]_\sim$ according to Definition \ref{big}\lref{sim} (b).
\end{rem}
%
%
%
%
%
%
%
%
%
\begin{Lemma}\label{semiglobal}
Let $r\ge 1$ and ${\cal F}$ be a graphical $C^r$ codim-1-foliation.
Then the following conditions are
equivalent:
\begin{enumerate}
\item \ilabel{loceul} $\delta ({\cal F})$ is locally eulerian.
\item \ilabel{semiglobalexistence} For each $x$$\in$$M$
there is some
neighbourhood $U$$\subset$$G$ of $p\bigl(\pi(x)\bigr)$, such that there is
$f$$\in$$C^r \bigl(\tilde U ,\real\bigr)$ with $df(x)\not=0$ in
$\tilde U$ and ${\cal F}\vert_{\tilde U}={\cal F}(df)$,
where $\tilde U :=\pi^{-1}\bigl(p^{-1}(U)\bigr)$.
\end{enumerate}
\end{Lemma}
Condition \iref{semiglobalexistence}
means the semiglobal existence of $f$.\\[0.5ex]
\begin{proof}
Let \iref{loceul} be satisfied. 
The transverse orientability of ${\cal F}$ and the existence of a piecewise $C^r$ function $f_0$ without extrema like described are rather technical to prove, but plausible. 
So we omit this part.
Furthermore it is clear that after any choice of a non-Hausdorff
leaf $F_0$ the function $f_0$ can be chosen $C^r$ around $F_0$ with $df_0\not=0$ on $F_0$. Now $C^r$-regularity and transverse orientability of ${\cal F}$ imply that there is a saturated neighbourhood
$U_0$ of $F_0$ and a function $\tilde f\in C^r_{\cal F}(M)$ such that $f_0|U_0=\tilde f|U_0$. 
For any $F_i\sim F_0$ holds automatically $d\tilde f\not=0$ on $F_i$: Let $T$ be a local transversal through $x\in F_0$ with injective projection to $M/{\cal F}$; this exists, since we assumed $E=\emptyset$. The function
$x\in F_t (t\in T)\mapsto \frac 1{d\tilde f(\frac{\partial}{\partial t})}$
would be in $C^{r-1}_{\cal F}(\sat(T))$
and not be extendable beyond $F_i$, if $d\tilde f|F_i=0$. So for all $F_i\sim F$ there is a saturated neighbourhood $U_i\supset F_i$ such that $d\tilde f\not= 0$ on $U_i$.
We set $U:=\bigcup_{F_i\in[F]}U_i$ and $f:=\tilde f|U$.

In order to prove the other direction, $\iref{semiglobalexistence}\imp\iref{loceul}$, the arguments are again technical, but straightforward.
%
\end{proof}
%
%
%
%
%
%
%
%
%
%
%
%
\begin{Theorem} \label{Hauptsatz} 
{\rm (Criterion for global existence) } Let ${\cal F}$ be a graphical $C^r$ codim-1-foliation.
Then the following
conditions are equivalent:
\begin{enumerate}
{\rm \item \hspace{0ex}}$\delta({\cal F})$ is globally eulerian.
{\rm \item \hspace{0ex}}There is some function $f\in C^r (M,\real)$ such that
 $df(x)\not=0$ in $M$ and ${\cal F}={\cal F}(df)$.
\end{enumerate}
\end{Theorem}
\begin{proof}
(i) Suppose $\delta({\cal F})$ to be globally eulerian.
Let ${\cal R}:={\cal V}\Gamma$ be partially ordered according to
\ref{bem} \iiref{drei}
and let ${\rm a}: {\cal V}\Gamma\to {\mathbb Q}$ be any mapping
satisfying
\[
\sigma < \sigma'\Longrightarrow {\rm a}(\sigma) < {\rm a}(\sigma')
\]
%
%
Moreover, for each $\sigma\in {\cal V}\Gamma$ choose
$\epsilon_\sigma\in (0,1)$
with
\begin{equation}\label{epsilonsigma}
\epsilon_\sigma < \frac{1}{3}\
{\rm min}\
\Big\{\ |{\rm a}(\sigma')-{\rm a}(\sigma)|\ \Big\vert\ \sigma'\in
{\cal V}\Gamma,\
\exists L\in {\cal E}\Gamma: \kappa(L)=\{\sigma,\sigma'\}\Big\}.
\end{equation}
Let $\tilde U_\sigma$ be an appropriate neighbourhood of $\sigma$, $x\in
s\in {\cal E}\gamma_\sigma$
and
$f_\sigma\in C^{r}\bigl(\tilde U_\sigma,\real\bigr)$ the corresponding mapping
as in condition (b) of Lemma \ref{semiglobal}.
Without any restriction we may
assume
$f_\sigma(x)\in [-\epsilon_\sigma, +\epsilon_\sigma]$ for each
$x\in\tilde U_\sigma$
and
$\tilde U_\sigma\cap\tilde U_{\sigma'}=\emptyset$ for
$\sigma'\in {\cal V}\Gamma\setminus\{\sigma\}$.
Now define
\[
\tilde f:\ \bigcup_{\sigma\in{\cal V}\Gamma} \tilde U_\sigma \to \real,\
\tilde f(x):= {\rm a}(\sigma)+f_\sigma(x)\, \mbox{ for }\, x\in\tilde U_\sigma.
\]
Let $L\in{\cal E}\Gamma$, $\kappa(L)=\{ \sigma,\sigma'\}$,
$L\in {\cal V}_+\gamma_\sigma\cap {\cal V}_-\gamma_{\sigma'}$. Thus, we
have ${\rm a}(\sigma) < {\rm a}(\sigma')$ and, using condition
(\ref{epsilonsigma}),
\[
\forall x\in\tilde U_\sigma,\
x'\in\tilde U_{\sigma'}:\quad \tilde f(x) <\tilde f(x').
\]
Moreover, let $u\in C^r\bigl( [0,1],M\bigr)$ be some transversal  curve
according
to \ref{Further notations} such that $\pi\circ u\in {\bf W}$ and
$L=[\pi\circ u]_\sim$ 
(For the definition of ${\bf W}$ see Construction \ref{Makrograph}).
For $x\in M$ with $x\in L$ according to (\ref{xinL}) we define
\[
f(x):=\phi\Bigl(\bigl(\pi\circ u\bigr)^{-1}\bigl(\pi(x)\bigr)\Bigr).
\]
for a suitably chosen strictly increasing function $\phi\in C^r([0,1])$.
After an appropriate reparametrization of $u$, we can achieve $f(x)=\tilde f_\sigma (x)$ for $x\in\tilde U_{\sigma}$, $\sigma\in {\cal V}\Gamma$.
Then $f$ has the desired properties.
\\[0.6ex]
(ii) Now suppose condition (b) to be fulfilled. 
From Lemma \ref{semiglobal} we know that
$\Gamma$ is locally eulerian. 
An ordering on ${\cal
V}\Gamma$ and therefore an orientation $\uparrow$ of $\Gamma$
according to Remark \ref{bem} (iii) is induced by $f$, given by
\[
\sigma <\sigma' \text{\quad if \quad } f(\xi) < f(\xi'),
\]
where 
$\sigma,\sigma'\in {\cal V}\Gamma$, $\xi,\xi'\in M$ with 
$p\bigl(\pi(\xi)\bigr)=\sigma$, $p\bigl(\pi(\xi')\bigr)=\sigma'$ (cf.\ 
Def.\ \ref{big}\lref{sim})\\
Furthermore, to see the bipartition of ${\cal V}\gamma_\sigma$, define ${\cal V}_+\gamma_\sigma$,
${\cal V}_-\gamma_\sigma$ as in Remark \ref{bem} (ii).
%
\end{proof}
%
%
%
%

%

%

%
%
%
%
%
\subsection{The case $b_1(M)=0$}

If ${\cal F}$ is a graphical foliation, then denote $\delta({\cal F})=(\Gamma,(\gamma_\sigma)_\sigma)$ the graphical configuration. Like above, we denote $\Gamma=({\cal V}\Gamma,
{\cal E}\Gamma)$ and similar for other graphs.
Let us define the main graph $\mu=\mu({\cal F})$ by
$${\cal V}\mu:=\bigcup_{\sigma\in{\cal V}\Gamma}{\cal V}\gamma_\sigma,\quad {\cal E}\mu:={\cal E}\Gamma\cup\bigcup_{\sigma\in{\cal V}\Gamma}{\cal E}\gamma_\sigma$$
and the corresponding boundary mapping
$$
\kappa_\mu:\,\ {\cal E}\mu\to S^2 {\cal V}\mu,\,\
\kappa_\mu(\Lambda):=\left\{\begin{array}{ll}\kappa_\sigma(\Lambda),\ &\mbox{ if }\Lambda\in{\cal E}\gamma_\sigma\\
\{\Lambda_\sigma,\Lambda_{\sigma'}\},\ &\mbox{ if }\Lambda\in{\cal E}\Gamma\end{array}\right. ,$$
where in the second case 
$\{\sigma,\sigma'\}=\kappa_\Gamma(\Lambda)=\kappa(\Lambda)$ 
(cf.\ (\ref{kappamacro}) in Construction \ref{Makrograph})
and $\Lambda_\sigma$ 
denotes $\Lambda$ as element of ${\cal V}\gamma_\sigma$, correspondingly for
$\sigma'$.

This means, that $\mu$ is the graph we obtain, if we replace the vertices of $\Gamma$ with their corresponding micrographs.

Let $\mu$ and $\Gamma$ also denote the topological space associated to the graph $\mu$ resp.\ $\Gamma$. Then it is clear that there is
a continuous map 
$$c:\mu\seq \Gamma,$$
which contracts the $\gamma_\sigma$ to the point $\sigma$ and is an isomorphism outside the micrographs.
Moreover, there is a continuous map
$$\tilde c:\mu\seq M/{\cal F},$$
which identifies the edges of $\gamma_\sigma$ with the corresponding non-Hausdorff leaves in $M/{\cal F}$, and is an isomorphism outside
the micrographs. This map satisfies $c=p\circ\tilde c$.
The fibres $c^{-1}(g)$ of $\tilde c$ are homeomorphic to intervals, hence 
they are simply connected. Outside of the micrographs $\tilde c$ is an
isomorphism and $M/{\cal F}$ is a one-dimensional 
manifold outside the non-Hausdorff points.
As usual, we denote $b_1(Y):=rk H^1(Y,{\bb Z})$ for a topological space $Y$.

By standard topological methods we obtain $H^1(M/{\cal F},{\bb Z})\cong H^1(\mu,{\bb Z})$. Applying this we arrive at

\begin{Lemma}\label{betti}If $\delta({\cal F})$ is the graphical configuration of a foliation, then
\begin{enumerate}
{\rm \item \ilabel{b1} $b_1(M/{\cal F})\le b_1(M)$,
\item \ilabel{b2} $b_1(\Gamma)+\sum_{\sigma\in{\cal V}\Gamma} b_1(\gamma_\sigma)=b_1(M/{\cal F}),$
\item \ilabel{b3} $\Gamma$} and all $\gamma_\sigma$ are simply connected, if $b_1(M/{\cal F})=0$.
\end{enumerate}
\end{Lemma}

Now we can derive the global existence theorem for manifolds with $b_1(M)=0$:

\begin{Theorem}\label{b1M=0}
Let $M$ be a manifold with $b_1(M)=0$ and ${\cal F}$ a graphical $C^r$ codim-1-foliation. Then there is a $C^r$ function $f:M\seq{\mathbb R}$ with
$df(x)\not=0$ everywhere and ${\cal F}={\cal F}(df)$.
\end{Theorem}

\begin{proof}By Lemma \ref{betti} we know, that $\Gamma$ and all $\gamma_\sigma$ are simply connected, so they contain no cycle. Since a graph is bipartite
if and only if every cycle has an even number of edges, we conclude, that all $\gamma_\sigma$ are bipartite, hence $\delta({\cal F})$ is locally eulerian.
Now let us define an orientation of $\Gamma$ inductively: choose any $\sigma\in
{\cal V}\Gamma$ and a bipartition
$${\cal V}\gamma_\sigma={\cal V}_+\gamma_\sigma\cup{\cal V}_-\gamma_\sigma\.$$
For any $L\in{\cal E}\Gamma^\sigma$ choose the orientation $\uparrow(L)$ 
and the bipartition of the neighboured micrographs so that
\begin{eqnarray*}\uparrow_1(L)=\sigma \wedge L\in{\cal V}_-\gamma_{\uparrow_2(L)}\ , &\mbox{ if }&L\in{\cal V}_+\gamma_\sigma,\\
\uparrow_2(L)=\sigma\wedge L\in{\cal V}_+\gamma_{\uparrow_1(L)}\ ,
&\mbox{ if }&L\in{\cal V}_-\gamma_\sigma.
\end{eqnarray*}
Now we are looking at the neighbouring micrographs and repeat this procedure.
This will give a complete system of orientation of $\Gamma$ and bipartitions
of the $\gamma_\sigma$, because $\Gamma$ is connected. 
We have to show that there are no inconsistencies. But if there are any,
say at $\gamma_{\sigma'}$, there would be two different
paths
\begin{eqnarray*}\sigma L_1\sigma_1...\sigma_{k-1}L_k\sigma'&\mbox{ resp.}&\sigma L_k\sigma_k...\sigma_{l-1}L_l\sigma'.
\end{eqnarray*}
Hence the path
$$\sigma L_1\sigma_1...\sigma_{k-1}L_k\sigma' L_l\sigma_{l-1}...\sigma_kL_k\sigma$$
would contain a circuit, what contradicts $b_1(\Gamma)=0$.

In particular, there can be no oriented circuits, so we conclude
that $\delta({\cal F})$ is globally eulerian. By Theorem \ref{Hauptsatz}
the proof is finished.
\end{proof}

\section{Applications}

\subsection{Tools}
%
This section is intended to provide tools for the investigation of
foliations. In particular, there is some simplification for the 
test of the graphicalness of a foliation.
\subsubsection{Closedness of leaves}

If we assume a leaf $F$ of a $C^\infty$ codim-1-foliation ${\cal F}$ on a manifold $M$ to be not closed (topological), then $F$ as a point in $M/{\cal F}$ is a non-Hausdorff
point. This we see by considering $\overline{F}$. It is easy to see that $\overline{F}$ is saturated, hence there exists a leaf $F'\subset\overline{F}\setminus F$.
Clearly, $F'$ and $F$ cannot be separated by open neighbourhoods. Hence Hausdorff leaves are closed.

If we want to verify the graphicalness of a foliation it turns out that
this gets much simpler, if on $M$ the Jordan Brouwer Separation Theorem holds for closed leaves,
that means:

\begin{defn}A manifold $M$ is called {\em foliated Jordan Brouwer separable} (fJBS), if 
for any $C^\infty$ codim-1-foliation ${\cal F}$ and closed leaf $F\in{\cal F}$ holds: $M\setminus F$ has at least two
connected components.
\end{defn}

\begin{Lemma}Any manifold $M$ with $H_1(M,{\bb Z})=0$ is fJBS.\end{Lemma}

\begin{proof}For any manifold $Z$ we denote by ${\cal A}_Z$ the sheaf of $C^\infty$ functions on $Z$.
Let $F$ be a closed leaf of a foliation ${\cal F}$ and $U(F)$ a small tubular neighbourhood. 
The condition $H_1(M,{\bb Z})$ implies that the part of the  Mayer-Vietoris sequence
$$0\seq H_0(U(F)\setminus F,{\bb Z})\seq H_0(M\setminus F,{\bb Z})\oplus H_0(F,{\bb Z})\seq H_0(M,{\bb Z})\seq 0$$
is exact. Since $H_0(U(F)\setminus F,{\bb Z})=H_0(|N_F|\setminus F,{\bb Z})$, if $|N_F|$ denotes the total space of the normal bundle of $F$, it is
enough to show $N_F\cong {\cal A}_F$: Then $|N_F|=F\times{\bb R}$, hence $|N_F|\setminus F=F\times({\bb R}\setminus\{0\})$.

The crucial point is now that $N_F$ extends to the line bundle $N_{\cal F}$ on $M$. Real line bundles are parametrised by
$H^1(M,{\bb Z}_2)$, however. By the universal coefficient lemma this group vanishes, so any line bundle on $M$ is trivial, in particular
$N_{\cal F}\cong {\cal A}_M$. This again implies $N_F\cong{\cal A}_F$. 

Note that we even proved that $M\setminus F$ has exactly two components.
\end{proof}

\begin{Lemma}\label{closed}Let $M$ be fJBS manifold and ${\cal F}$ a codim-1-foliation on $M$. Then the set
$$C:=\left\{F\in \left. M/{\cal F}\ \right|\ F \mbox{ is a closed leaf}\right\}$$
is closed.
\end{Lemma}

The proof is straightforward and hence left to the reader.

Note that a leaf $F\subset M$ is closed if and only if $F\in M/{\cal F}$
is a closed point.

The lemma implies immediately:

\begin{cor}Let $M$ be a fJBS manifold. Then a regular $C^r$ 
codim-1-foliation ${\cal F}$ is graphical if and only if 
$p(S)$ is a discrete set, and for every $g\in G$ the fibre $p^{-1}(g)$ is a finite set, where $S$
is like in Definition \ref{big}\lref{sim} the set of non-Hausdorff points.
\end{cor}

\begin{proof} By the discreteness assumption, $G\setminus p(S)$ is open and dense, hence $H:=M/{\cal F}\setminus S$ is open, since $p$ is continuous, and dense,
since the fibres are finite.
But $H$ is the set of Hausdorff points, hence $H\subset C$. This means by Lemma \ref{closed} that $C=\overline C=M/{\cal F}$, hence all leaves are closed.
\end{proof}
 
%
%
\subsubsection{Characterization of infinitely near leaves}
Here we give criteria for constructing the graphical configuration directly from the one-form. We suppose $M$ to be Riemannian. Then we have the notion
of the normal bundle $N_{\cal F}$ to ${\cal F}$, defined by
$$N_{{\cal F},x}:=(T_xF_x)^\perp,$$
where $F_x$ is the leaf containing $x$. A normal curve $b:[0,1]\seq M$ then is a curve with $\dot{b}(t)\in N_{{\cal F},b(t)}$ for all $t\in[0,1]$, i.e.
an integral curve of the normal bundle.
\par
\begin{figure}[h!]
\begin{minipage}[h]{7.0cm}
\begin{center}
\includegraphics[height=13ex]{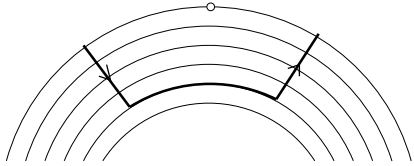}
\caption{A normal path}
\end{center}
\end{minipage}
\begin{minipage}[h]{7.0cm}
Let us define a {\it normal path} $\beta$ connecting $x\in F$ and $x'\in F'$ as the sum $\nu+b+\nu'$, where
$\nu,\nu':I\seq M$ are paths whose image is contained in the curve normal to ${\cal F}$ meeting $x$ resp $x'$ and
$\nu(0)=x,\nu'(1)=x'$, and $b$ is a path whose image is contained in a leaf $\tilde F$ and $\nu(1)=b(0), b(1)=\nu'(0)\.$  
\end{minipage}
\end{figure}
\par
\begin{Lemma}\label{int}Let $\omega$ be an integrable one-form without zeroes
on the Riemannian manifold $(M,g)$ and ${\cal F}={\cal F}(\omega)$ the induced foliation. Then the classes of any two leaves 
$F$ and $F'$ are infinitely near in $M/{\cal F}$ if and only if
for some (and hence all) $x\in F, x'\in F'$ and all $\varepsilon>0$ exists a 
piecewise continuous differentiable normal path $\beta:[0;1]\seq M$ with 
$\beta(0)=x, \beta(1)=x'$ and
$$\int_\beta |\omega|:=\int_0^1|\beta^*(\omega)(\frac{d}{ds})|ds<\varepsilon\.$$
\end{Lemma}

\begin{proof}$'\Rightarrow'$: If $F,F'$ are infinitely near, then by definition, any saturated neighbourhoods 
$$U_\delta(F):=\left\{\tilde F\in M/{\cal F}\ \left| \inf_{x\in F,\tilde x\in\tilde F}(dist_g(x,\tilde x))\right.<\delta\right\}$$ and, analogously defined,
$U_\delta(F')$ satisfy
$$U_\delta(F)\cap U_\delta(F')\not=\emptyset\.$$
So, if 
$B_\delta(x)$ is the geodesic ball of radius $\delta$ around $x$ and 
$l_x:[0;1]\seq B_{\delta}(x)$
is a curve normal to ${\cal F}$ with $l_x(0)=x$, 
and analogously for $l_{x'}$, such that
there is a leaf $\tilde F$ with $l_x(1)\in\tilde F\ni l_{x'}(1)$, then let $\alpha:[0;1]\seq\tilde F$ be a path from $l_x(1)$ to $l_{x'}(1)$ and finally
$$\beta:= l_x+\alpha+l_{x'}\.$$
Clearly, $\beta$ is a normal path. If we choose $\delta$ sufficiently small, we can achieve
\begin{equation}\label{eps}\int_\beta |\omega|<\varepsilon.\end{equation}


$'\Leftarrow'$: Let $U$ and $U'$ be small saturated neighbourhoods of $F$ resp.\ $F'$. 
It is an easy argument, that for sufficiently small $\varepsilon>0$, the normal paths $\beta_\varepsilon$ are contained
in $U\cap U'$. In particular, $U\cap U'\not=\emptyset$.

\end{proof}

\begin{Lemma}\label{nu}
Let $\omega$ be an integrable one-form without zeroes on the complete Riemannian manifold $M$ and ${\cal F}={\cal F}(\omega)$ 
the induced foliation. $F$ and $F'$ are infinitely near in $M/{\cal F}$ if and only if for any
saturated neighbourhoods $U(F), U'(F')$ there is a by arc length parametrized
curve $\nu:{\bb R}^+\seq M$
which is normal to the foliation, with $\nu({\bb R}^+)\cap F=\emptyset,$ $\nu({\bb R}^+)\cap F'=\emptyset$, but
$\nu({\bb R}^+)\cap U(F)\cap U'(F')\not=\emptyset$.
\end{Lemma}

\begin{proof}$'\Leftarrow'$: If $\nu({\bb R}^+)\cap U(F)\cap U'(F')\not=\emptyset$ for saturated neighbourhoods $U(F)$, $U'(F')$,
then $U(F)\cap U'(F')\not=\emptyset$, hence $F$ and $F'$ are infinitely near.

$'\Rightarrow'$: For any neighbourhoods $U(F), U'(F')$ there is an $\varepsilon>0$ such that
$$U\cap U'\subset U(F)\cap U'(F'),$$ 
with saturated sets $U:=\pi^{-1}\circ\pi\circ\nu_\varepsilon([0,1]), U':=\pi^{-1}\circ\pi\circ\nu'_\varepsilon([0,1])$ containing $F$ resp.\ $F'$ 
where $\beta_\varepsilon=\nu_\varepsilon+b+\nu'_\varepsilon$ denotes a normal path satisfying (\ref{eps}) like in the proof of Lemma \ref{int}. 
Furthermore fix a leaf $\tilde F\subset U\cap U'$.
Now look at the continuous map induced by the transport along the normal curves to ${\cal F}$
$$n:H\seq \tilde F,$$
where $H$ 
is the locus of points $p\in F\cup F'$ such that the normal curve through $p$ meets $\tilde F$.
By construction, $H\cap F\not=\emptyset$, $H\cap F'\not=\emptyset$, so $H$ has at least two connected components. 
Since through a point passes exactly one normal curve,
$n$ is injective, hence $n(H)$ has at least two connected components. Thus $n(H)\not=\tilde F$. Now take a point $x\in\tilde F\setminus
n(H)$ and look at the normal curve $\nu:{\bb R}^+\seq M$ with $\nu(0)=x$. By completeness of $M$ it can be achieved, that
$\nu$ is parametrized by arc length with domain ${\bb R}^+$. This curve satisfies every condition of the lemma.
\end{proof}
\begin{figure}[ht]
\begin{center}
\includegraphics[height=15ex]{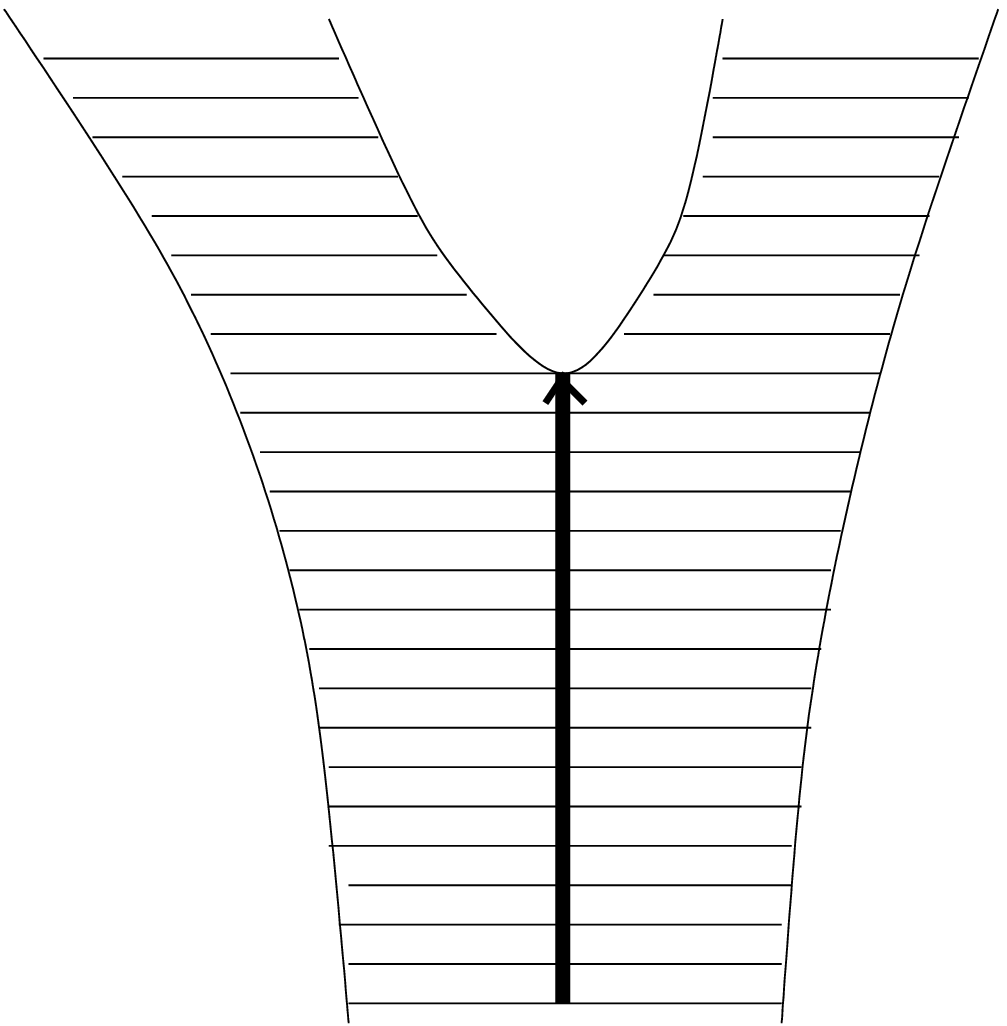}\hspace{4ex}
\includegraphics[height=15ex]{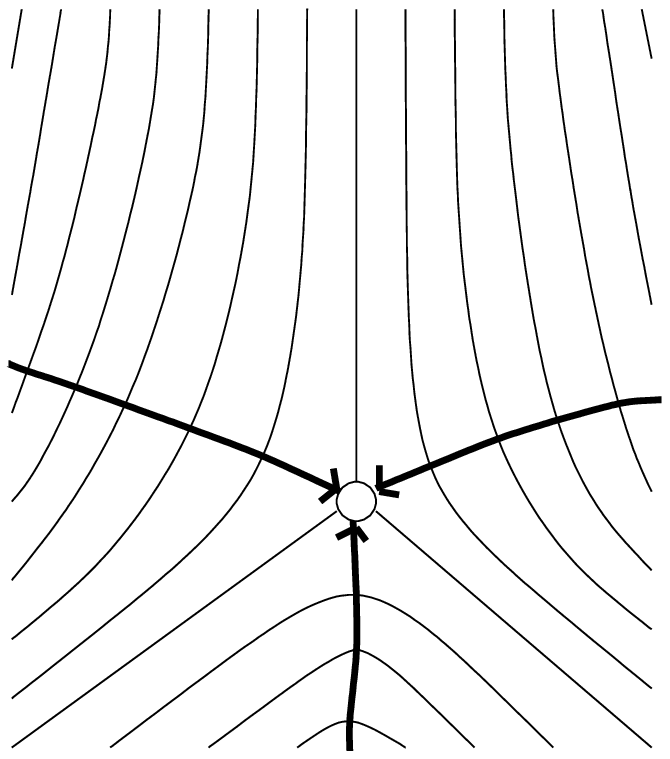}\hspace{4ex}
\includegraphics[height=10ex]{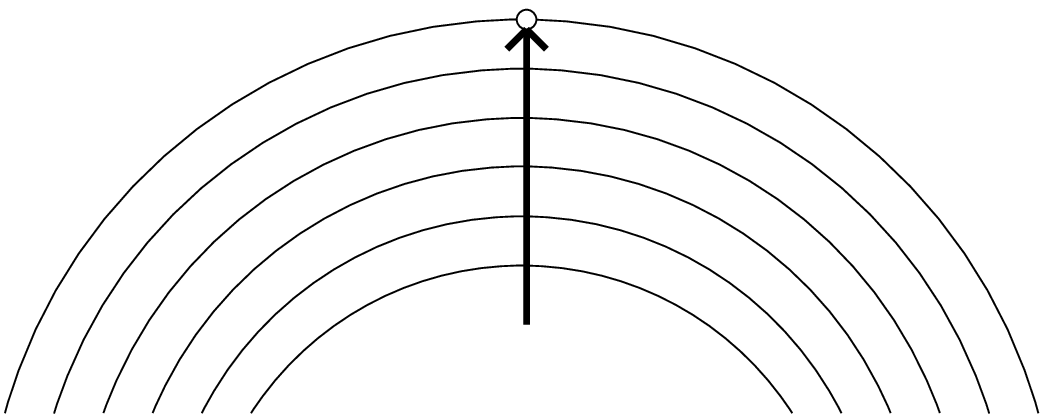}
\end{center}
\caption{Examples for the curve $\nu$ in Lemma \ref{nu}}
\end{figure}

% ...
\subsection{Special cases}
%
In this section we present special cases, in which the decision of the existence of an Eulerian multiplier is easier than in other cases. These examples are
compact leaves of the foliation and a cylinder form of the manifold.
%
\subsubsection{The case of compact leaves}

Let us discuss the case that the leaves are compact. It will turn out that under compactness of the leaves $\hol(F)=1$ is equivalent to
the existence of an Eulerian multiplier.

\begin{Theorem}\label{cp}Let $M$ be a non-compact $C^r$ manifold and ${\cal F}$ a codimension one $C^r$ foliation, such that all leaves $F$ are compact 
and $\hol(F)=1$. Then there is a $C^r$ function $f:M\seq{\real}$ satisfying ${\cal F}={\cal F}(df)$ and $df(x)\not=0$ in $M$.
\end{Theorem}

\begin{proof} By \cite[Ch.\ IV, Lemma 6]{ccaln}, for any leaf $F$ there is a saturated neighbourhood
$U(F)$ together with a commutative diagram
$$\xymatrix{
U(F)\ar[d]^\pi\ar[r]^h & I\times F\ar[d]^{pr}\\
V\ar[r]^{i_V} & I
},$$
where $V$ is an open set in $M/{\cal F}$, $I$ is an interval in ${\real}$, $h$ is a $C^r$ diffeomorphism and $i_V$ a
homeomorphism. This suggests a natural Hausdorff manifold structure on $M/{\cal F}$. By the non-compactness of $M$ we deduce that $M/{\cal F}$ is not compact and hence
diffeomorphic to ${\bb R}$. Thus the projection is the desired function $f$.
\end{proof}

This theorem allows a stronger result, if $\pi_1(M)$ is finite, much similar to Corollary \ref{b1M=0}. Indeed, the proof of Theorem \ref{cp} used in the following shows with some extra
arguments presented below, that
${\cal F}$ is a graphical foliation, if $\pi_1(M)$ is finite and all leaves are compact:

\begin{cor}\label{pi1=0}Let $M$ be a non-compact Riemannian manifold with finite fundamental group $\pi_1(M)$. If ${\cal F}$ is a $C^r$ codimension one foliation with only compact leaves, then
there is a $C^r$ function $f:M\seq \real$ with $df(x)\not= 0$ for all $x\in M$ such that ${\cal F}={\cal F}(df)$.
\end{cor}

\begin{proof}(i) As a first case, let us assume $\pi_1(M)=0$.

By Lemma \ref{hol} we know that $\hol(F)\in\{1,{\bb Z}_2\}$ for all leaves $F$.
Now we can follow the proof of the Corollary in \cite[Ch.\ IV, \S 5]{ccaln} to see that $\hol(F)=1$ if the normal bundle $N_{\cal F}$ to ${\cal F}$ 
is orientable. But since $\pi_1(M)=0$ every
vector bundle over $M$ is orientable; in particular, so is $N_{\cal F}$. 
So we now can apply Theorem \ref{cp} to complete the proof.

(ii) If $|\pi_1(M)|=d$ is finite, we look at the universal cover $\tilde M\stackrel{u}{\seq}M$, what is a $d:1$-cover. The foliation
$\tilde{\cal F}:=u^*{\cal F}$ then has also only compact leaves, since for a leaf $\tilde F\in\tilde{\cal F}$ the map $u|\tilde F:\tilde F\seq F$,
with $F\in{\cal F}$ is a finite cover also (maybe of lower order). 
So we can apply the first part for $(\tilde M,\tilde{\cal F})$ and see that
$\tilde{\cal F}$ is graphical. This implies immediately that ${\cal F}$ is of finite type. In order to see that ${\cal F}$ is regularly $C^r$, by definition we may assume
that ${\cal F}$ is transversely orientable, hence given by an integrable one-form $\omega$. By the first part, $u^*\omega=\tilde fd\tilde g$ with
nonvanishing $\tilde f$; so we may assume $\tilde f>0$. Applying standard covering theory
yields $\omega=fdg$ with $$f\circ u=\frac 1d(\sum_{\phi\in Deck(\tilde M,M)}\frac{1}{\tilde f\circ\phi})^{-1},\,\, g\circ u=\frac 1d\sum_{\phi\in Deck(\tilde M,M)}\tilde g\circ\phi.$$
In particular, ${\cal F}$ is given by a closed one-form, so it is regularly $C^r$.
Now we know that ${\cal F}$ is graphical. The finiteness of $\pi_1(M)$ implies $b_1(M)=0$.
Now we apply Theorem \ref{b1M=0} and thereby proved the claim.  
\end{proof}

%
\begin{ex}{}If we take $M={\bb R}^n\setminus\{0\}$, $\omega$ a one-form on $M$ such that all leaves of ${\cal F}(\omega)$ are diffeomorphic to $S^{n-1}$, then the existence of 
an Eulerian multiplier is given by Corollary \ref{pi1=0}, if $n\ge 3$.
For $n=2$ the statement can be proved by elementary means.
\end{ex}

%
\subsubsection{ The case of an infinite cylinder}
Finally, we mention an obvious result.
\begin{Theorem}
Let $D\subset\real^n$ be a simply connected domain, $M:=D\times\real$. Moreover, let
$v=(v^1,...,v^{n+1})\in C^r(M,\real^{n+1})$ such that $v^{n+1}(x)\not=0$ for all $x\in M$ and $|\frac{v^{i}}{v^{n+1}}|\le C$ for all $i$ and some constant $C$. 
Let $\omega=:\sum_{k=1}^{n+1} v^k dx_k$ be integrable. 
There are functions $f\in C^r(M,\real)$ and 
$\lambda\in C^{r-1}(M,\real)$ satisfying $\lambda(x)\not=0$ in $M$ and
\[
df=\lambda \omega.
\]
\end{Theorem}

\begin{proof}The assumptions imply that for every leaf the projection $F\seq D$ is onto and locally an isomorphism (in the abstract manifold structure of $F$). 
Since $\pi_1(D)=0$, the projection
has to be a global isomorphism. In particular, $M/{\cal F}\cong{\bb R}$ and the claim follows.
\end{proof} 

{\bf Remark:} Note that, of course, the property $\omega\land d\omega=0$ and
the existence of Euler's multiplier are invariant under $C^r$- diffeomorphisms
$\Phi: M\to M'$, $r\ge 2$. Thus, the result of this theorem  still holds
if $M$ and the conditions are ``deformed'' consistently.
%
%
%
%
%
%
%
%
%
%

%

%
\end{document}